%% file: sequential-multiple-testing.tex
\documentclass[twoside,11pt]{article}

\usepackage{graphicx, subfigure}
\usepackage[top=1in,bottom=1in,left=1in,right=1in]{geometry}
\usepackage[sort&compress]{natbib} 
\usepackage{amsfonts, amsmath, amssymb, amsthm, constants, bbm}
\usepackage{MnSymbol}
\usepackage{mathtools}
\usepackage[mathscr]{eucal}

\usepackage{color}
\definecolor{darkred}{RGB}{100,0,0}
\definecolor{darkgreen}{RGB}{0,100,0}
\definecolor{darkblue}{RGB}{0,0,150}

\usepackage{hyperref}
\hypersetup{colorlinks=true, linkcolor=darkred, citecolor=darkgreen, urlcolor=darkblue}
\usepackage{url}

\input notation.tex

\definecolor{purple}{rgb}{0.4,.1,.9}
\definecolor{new}{rgb}{0.5,0.1,0.1}

\newcommand\blfootnote[1]{%
  \begingroup
  \renewcommand\thefootnote{}\footnote{#1}%
  \addtocounter{footnote}{-1}%
  \endgroup
}

\begin{document}
\thispagestyle{empty}

\title{Sequential Multiple Testing}
\author{Shiyun Chen \and Ery Arias-Castro}
\date{University of California, San Diego}
\maketitle

\begin{abstract}
We study an online multiple testing problem where the hypotheses arrive sequentially in a stream.  The test statistics are independent and assumed to have the same distribution under their respective null hypotheses. We investigate two procedures LORD and LOND, proposed by \citep{javanmard2015online}, which are proved to control the FDR in an online manner.  In some (static) model, we show that LORD is optimal in some asymptotic sense, in particular as powerful as the (static) Benjamini-Hochberg procedure to first asymptotic order.  We also quantify the performance of LOND.  Some numerical experiments complement our theory.  
\end{abstract}

\blfootnote{Both authors are with the Department of Mathematics, University of California, San Diego, USA.  Contact information is available \href{http://www.math.ucsd.edu/people/graduate-students/}{here} and \href{http://math.ucsd.edu/\~eariasca}{here}.  This work was partially supported by grants from the US National Science Foundation (DMS 1223137).}


\section{Introduction} \label{sec:intro}
Multiple testing is now a well-established area in statistics and arises in almost every scientific field \citep{dudoit2007multiple, dickhaus2014simultaneous, roquain2011type}.
In this paper, we consider a scenario where infinitely many hypotheses $\cH = (\H_1, \H_2, \H_3, \dots)$ arrive sequentially in a stream with corresponding P-values $P_1,P_2, P_3, \dots$, and we are required to decide whether we accept or reject $\H_i$ only based on $P_1, \dots, P_i$. We propose to use the recent sequential multiple testing procedures of \citep{javanmard2015online} which control the FDR in an online manner, and study the asymptotic optimality properties of these methods in the context of sparse mixture asymptotically generalized Gaussian model (see \defref{AGG}) which the normal model often used as benchmark in various works on multiple testing.

\subsection{The risk of a multiple testing procedure}
Consider a setting where we want to test an ordered infinite sequence of  null hypotheses, denoted $\cH = (\H_1, \H_2, \H_3, \dots, $), where at each step $i$ we have to decide whether to reject $\H_i$ having access to only previous decisions.  The test that we use for $\H_i$ rejects for large positive values of a statistic $X_i$.  Throughout, we assume that test statistics are all independent.  Denote the collection of the first $n$ hypotheses in the stream by $\cH(n) = (\H_1, \dots, \H_n)$, and the vector of first $n$ test statistics by $\bX(n) = (X_1, \dots, X_n)$. 
Let $\Phi_i$ denote the survival function\footnote{In this paper, the survival function of a random variable $Y$ is defined as $y \mapsto \P(Y \ge y)$.} of $X_i$ and $\boldsymbol \Phi(n) = (\Phi_1, \dots, \Phi_n)$.
We assume that the corresponding P-values can be computed (or at least approximated).  The simplest such case is when $\H_i$ is a singleton, $\H_i = \{\Phi_i^{\rm null}\}$, and the null distributions $\Phi_1^{\rm null}, \Phi_2^{\rm null}, \dots, $ are known.  In that case, the $i$-th P-value is defined as $P_i = \Phi_i^{\rm null}(X_i)$, which is the probability of exceeding the observed value of the statistic under its null distribution.  

Let $\F$ index all the false null hypotheses in the stream, and let $\F_n \subset [n] := \{1, \dots, n\}$ index the false null hypotheses in the first $n$ hypotheses, meaning
\beq
\F_n = \{i \in [n] : \Phi_i \notin \H_i\}.
\eeq
A multiple testing procedure $\R$ takes the infinite sequence of test statistics $\bX$ and returns a subset of indices representing the null hypotheses that the procedure rejects.  Given such a procedure $\R$, the false discovery rate is defined as the expected value of the false discovery proportion \citep{benjamini1995controlling}
\beq
\fdr_n(\R) = \E_{\boldsymbol\Phi}[\fdp_n(\R(\bX))], \quad \fdp_n(\R) := \frac{|\R(\bX(n)) \setminus \F_n|}{|\R(\bX(n))|},
\eeq
where we denoted the cardinality of a set $\A \subset [n]$ by $|\A|$ and use the convention that $0/0 = 0$.
While the FDR of a multiple testing procedure is analogous to the level or size of a test procedure, the false non-discovery rate (FNR) plays the role of power and is defined as the expected value of the false non-discovery proportion\footnote{This definition is different from that of \cite{genovese2002operating}.}
\beq
\fnr_n(\R) = \E_{\boldsymbol\Phi}[\fnp_n(\R(\bX))], \quad \fnp_n(\R) := \frac{|\F_n \setminus \R(\bX(n))|}{|\F_n|}.
\eeq
In analogy with the risk of a test --- which is defined as the sum of the probabilities of type I and type II error --- we define the risk of a multiple testing procedure $\R$ as the sum of the false discovery rate and the false non-discovery rate
\beq\label{risk}
{\rm risk}_n(\R) = \fdr_n(\R) + \fnr_n(\R).
\eeq

\begin{rem}
The procedure that never rejects and the one that always reject both achieve a risk of 1, so that any method that has a risk exceeding 1 is useless.
\end{rem}

\subsection{More related work}
The literature on multiple testing is by now vast.  Only more recently, though, have multiple testing procedures been proposed for the sequential setting. In the context of testing $J>2$ (fixed) null hypotheses about $J$ sequences of data streams of arbitrary size, \citep{bartroff2014multiple} proposes general stepup and stepdown procedures which provide control of the simultaneous generalized type I and II error rates. See also \citep{bartroff2014sequential} for procedures controlling the type I and II FWER's, and \citep{bartroff2013sequential} for procedures controlling the FDR and FNR (defined differently). 

Another situation also considered in literature is where the hypotheses are ordered based on prior information on how promising each hypothesis is. In this context, \citep{g2016sequential} develops two rules (FowardStop and StrongStop) to choose the number of hypotheses to reject which are shown to control the FDR.  A variation of StrongStop rule can also be applied in sequential model selection in regression model. \citep{barber2015controlling} proposes the Sequential stepup procedure (SeqStep) which also guarantees FDR control under independence. \citep{li2016accumulation} develops a broader class of ordered hypotheses testing procedures under such setting, called \textit{accumulation tests}, which generalize the existing two methods (FowardStop and SeqStep). 
\citep{lei2016power} derives an improved version of Selective SeqStep, called Adaptive SeqStep. 
See \citep{fithian2014optimal, fithian2015selective, lockhart2014significance} for more methods and applications in selective sequential model selection.

Still in the sequential setting, \citep{foster2008alpha} develops an alpha-investing procedure which provides uniform control of mFDR (a weaker control than FDR control) in online testing under some condition. The alpha-investing rule spends some of the wealth to perform each test and earns more wealth each time a discovery is made.  
\citep{aharoni2014generalized} provides a broader class of online procedures called generalized alpha-investing and also establish mFDR control. \citep{javanmard2015online} proposes two procedures called LOND and LORD algorithms which control both FDR and mFDR in online testing. We refer to \secref{LORD} and \ref{sec:LOND} for more details of rules and discuss their asymptotic risk in our context. More generally, \citep{javanmard2016online} studies generalized alpha-investing rules and obtains conditions for FDR control under a general dependence structure of test statistics. They also develop modified procedures for online control of the false discovery exceedance.

In the present paper we study some asymptotic power properties of the LORD and LOND methods, complement the results of \citep{javanmard2015online}.  
This paper is a continuation of our previous work in the static setting \citep{ariaschen2016distribution}, where an asymptotic oracle risk bound for multiple testing is obtained, and both the method of \cite{benjamini1995controlling} and the distribution-free method of \cite{barber2015controlling} are proved to achieve that bound.
Various other oracle bounds and corresponding optimality results for multiple testing procedures are available in the literature; see, for example, \citep{genovese2002operating, 
	sun2007oracle,
	storey2007optimal,
	bogdan2011asymptotic,
	neuvial2012false, 
	meinshausen2011asymptotic,
	ji2012ups,
	jin2014rare,
	butucea2015variable}.

\subsection{Content}
The rest of the paper is organized as follows.
In \secref{normal-model} we consider the normal location model and derive the performance of LORD under this model. Generalizing this model, in \secref{AGG model} we consider a nonparametric Asymptotic Generalized Gaussian model. We analyze the asymptotic performance of the LORD and LOND procedures of \cite{javanmard2015online} under this model in \secref{LORD} and \secref{LOND}. We present some numerical experiments in \secref{numerics}. All proofs are gathered in \secref{proofs}.

\section{Methods} \label{sec:methods}

We describe the LORD and LOND procedures of \cite{javanmard2015online}, which are the methods we study in this paper.
Recall that $\H_1, \H_2, \dots$ are tested sequentially and that $P_i$ denotes the P-value corresponding to the test of $\H_i$.   These two procedures, and most others, work as follows: set a significance level $\alpha_i$ based on $P_1, \dots, P_{i-1}$ (except for $\alpha_1$ which is set beforehand) and reject $\H_i$ if $P_i \le \alpha_i$.  The LORD and LOND methods vary in how they set these thresholds, although they both start with a sequence of the form
\beq\label{lambda}
\lambda_i \ge 0 \text{ such that } \sum_{i = 1}^{\infty} \lambda_i = q,
\eeq
where $q$ denotes the desired FDR control level
In what follows, we stay close to the notation used in \citep{javanmard2015online}.

\subsection{The LORD method}
Based on a chosen sequence \eqref{lambda}, the LORD algorithm --- which stands for (significance) Levels based On Recent Discovery --- sets the sequential significance levels $(\alpha_i)_{i = 1}^{\infty}$ as follows:
\beq \label{LORD}
\alpha_i = \lambda_{i - t_i}, \quad t_i = \max\{l < i: \H_l \text{ is rejected}\}, 
\eeq
with $t_1 := 0$.

In \citep{javanmard2016online} the LORD algorithm is shown to control FDR at a level less than or equal to $q$ in an online fashion, specifically,
\beq\label{online-fdr}
\sup_{n \ge 1} \fdr_n(\R) \le q,
\eeq 
if the P-values are independent. 
More generally, \cite{javanmard2016online} study a class of monotone generalized alpha-investing procedures (which includes LORD as a special case) and prove that any rule in this class controls the cumulative FDR at each stage  provided the P-values corresponding to true nulls are independent from the other P-values.  

\subsection{LOND}
Based on a chosen sequence \eqref{lambda}, the LOND algorithm --- which stands for (significance) Levels based On Number of Discovery --- sets the sequential significance levels $(\alpha_i)_{i = 1}^{\infty}$ as follows:
\beq \label{LOND}
\alpha_i = \lambda_i (D(i-1) + 1).
\eeq
where $D(n)$ denotes the number of discoveries in $\cH(n) = (\H_1, \dots, \H_n)$, with $D(0) := 0$.

In \citep{javanmard2015online} the LOND is shown to control FDR at level less than or equal to $q$ everywhere in an online manner, the same as \eqref{online-fdr}, if the P-values are independent.

\section{Models} \label{sec:models}

In this paper we study the FNR of each of the LORD and LOND methods of \cite{javanmard2015online} on the first $n$ hypotheses as $n \to \infty$.  
As benchmark, we use the oracle that we considered previously\citep{ariaschen2016distribution} for the static setting defined by these $n$ hypothesis testing problems.  For the reader not familiar with that paper, at least in the models that we consider, this turns out to be asymptotically equivalent to applying the Benjamini-Hochberg (BH) method to the first $n$ hypotheses.  Note that the latter accesses all the first $n$ hypotheses at once and is thus not constrained to be sequential in nature.

The static setting we consider is that of a location mixture model.  
We assume that we know the null distribution function $\Phi$, assumed to be continuous for simplicity. We then assume that the test statistics are independent with respective distribution $X_i \sim \Phi_i = \Phi(\cdot - \mu_i)$, where $\mu_i = 0$ under the null $\H_i$ and $\mu_i > 0$ otherwise.
Both minimax and Bayesian considerations lead one to consider a prior on the $\mu_i$'s where a fraction $\eps$ of the $\mu_i$'s are  randomly picked and set equal to some $\mu > 0$, while the others are set to 0.  The prior is therefore defined based on $\eps$ and $\mu$, which together control the signal strength.  The P-value corresponding $\H_i$ is $P_i := \bar\Phi(X_i)$, where $\bar\Phi := 1 - \Phi$ is the null survival function.

\subsection{The normal model} \label{sec:normal-model}

As an emblematic example of the distributional models that we consider in this paper, let $\Phi$ denote the standard normal distribution.  Assume as above that $X_i \sim \Phi$ under $\H_i$ and $X_i \sim \Phi(\cdot - \mu)$ otherwise.  Thus, under the each null hypothesis, the corresponding test statistic is standard normal, while that statistic is normal with mean $\mu$ and unit variance otherwise.
This is the model we consider in \citep{ariaschen2016distribution} and the inspiration comes from a line of research on testing the global null $\bigcap_i \H_i$ in the static setting \citep{ingster1997some,IngsterBook,donoho2004higher}.
As in this line of work, we use the parameterization pioneered by \cite{ingster1997some}, namely
\beq \label{parameterization}
\eps = n^{-\beta} \text{ and }
\mu = \sqrt{2r \log n}.
\eeq

In the static setting, we know from our previous work \citep{ariaschen2016distribution} that any threshold-type procedure has risk tending to~1 as $n \to \infty$ when $r < \beta$ are fixed.  We also know that the BH method with FDR control at $q \to 0$ slowly has risk tending to 0 when $r > \beta$ are fixed.  In fact, these results are derived in the wider context of an asymptotically generalized Gaussian model, which we consider later.
Thus $r = \beta$ is the static selection boundary.


\begin{rem}
\citep{javanmard2015online} compared the power of their procedures in terms of lower bounds on the total discovery rate under the same mixture model but with a fixed mixture weight $\eps$.  In contrast, here we focus on the setting where $\eps \to 0$, meaning that the fraction of false null hypotheses (i.e., true discoveries) is negligible compared to the total number of null hypotheses being tested.
\end{rem}

\subsection{Asymptotically generalized Gaussian model}
\label{sec:AGG model}

Beyond the normal model, we follow \citep{ariaschen2016distribution, donoho2004higher} and consider other location models where the base distribution has a polynomial right tail in log scale.  

\begin{Def} \label{def:AGG}
A survival function $\bar\Phi = 1-\Phi$ is asymptotically generalized Gaussian (AGG) on the right with exponent $\gamma > 0$ if $\lim_{x \to \infty} x^{-\gamma} \log \bar{\Phi}(x) = -1/\gamma$.
\end{Def}

The AGG class of distributions is nonparametric and quite general.  
It includes the parametric class of generalized Gaussian (GG) distributions with densities $\{\psi_\gamma, \gamma > 0\}$ given by $\log \psi_\gamma(x) \propto -|x|^\gamma/\gamma$, which comprises the normal distribution ($\gamma = 2$) and the double exponential distribution ($\gamma=1$).  We assume that $\gamma \ge 1$ so that the null distribution has indeed a sub-exponential right tail.

\begin{rem}
We note that the scale (e.g., standard deviation) is fixed, but this is really without loss of generality as both the LORD and LOND methods are scale invariant.  This is because the P-values are scale invariant.
\end{rem}

The model is the same at that considered in \secref{normal-model} except that $\Phi$ is an unspecified (but known to the statistician) AGG distribution with parameter $\gamma \ge 1$. 
As in \citep{donoho2004higher}, we use the following parameterization
\beq \label{n-mu}
\eps = n^{-\beta} \text{ and }
\mu = \left( \gamma r \log n  \right)^{1/\gamma},
\eeq
where $r \ge 0$ and $\beta \in (0,1)$ are always assumed fixed.

\section{Performance analysis}

In this section we analyze the performance of the LORD and LOND methods in the static setting described earlier.  
Recall that $q$ denotes the desired FDR control level.  Typically $q$ is set to a small number, like $q = 0.10$.  In this paper we allow $q \to 0$ as $\eps \to 0$, but slowly.  Specifically, we always assume that
\beq\label{q}
\text{$q = q(n) > 0$ and $n^a q(n) \to \infty$ for all fixed $a > 0$.}
\eeq

\subsection{The performance of LORD} \label{sec:LORD}

We first establish a performance bound for LORD.
It happens that, despite required to control the FDR in an online fashion, LORD achieves the static selection boundary when desired FDR control is appropriately set.

\begin{thm}[Performance bound for LORD] \label{thm:AGG-LORD}
Consider a static AGG mixture model with exponent $\gamma \ge 1$ parameterized as in \eqref{n-mu}.  Assume that we apply LORD with $(\lambda_i)_{i = 1}^{\infty}$ defined as $\lambda_i \propto i^{-\nu}$ with $\sum_{i = 1}^{\infty} \lambda_i = q$, where $\nu > 1$ and $q$ satisfies \eqref{q}. 
If $r > \nu \beta$, the LORD procedure has $\fnr_n \to 0$ as $n \to \infty$.  
In particular, if $q \to 0$, then it has risk tending to 0.
\end{thm}

Note that the latter part comes from the fact that the LORD procedure controls of the FDR at the desired level $q$ as established in \citep{javanmard2015online} in the more demanding online setting.
In essence, therefore, LORD (with a proper choice of $\nu$ above) achieves the static oracle selection boundary $r = \beta$.

\begin{rem}\label{rem:lambda}
Assume that, instead, we apply LORD with any decreasing sequence $(\lambda_i)_{i = 1}^{\infty}$ satisfying $\sum_{i = 1}^{\infty} \lambda_i = q$ and
\beq \label{newlambda}
i^\nu \lambda_i \to \infty , \text{ for any fixed } \nu > 1.
\eeq
Then the conclusions of \thmref{AGG-LORD} remain valid.  In particular, such a choice of sequence (e.g., $\lambda_i \propto (\log i)^2/i$) adapts to the (usually unknown) values of $r$ and $\beta$.  (We provide details in \secref{proofs}.)
\end{rem}

\subsection{The performance of LOND} \label{sec:LOND}

We now turn to LOND and establish a performance bound under the same setting.  

\begin{thm} \label{thm:AGG-LOND}
Consider a static AGG mixture model with exponent $\gamma \ge 1$ parameterized as in \eqref{n-mu}.  Assume that we apply LOND with $(\lambda_i)_{i = 1}^{\infty}$ defined as $\lambda_i \propto i^{-\nu}$ with $\sum_{i = 1}^{\infty} \lambda_i = q$, where $\nu > 1$ and $q$ satisfies \eqref{q}. 
If $r > \beta + (\nu^{1/\gamma} - r^{1/\gamma})^\gamma + \nu -1$, the LORD procedure has $\fnr_n \to 0$ as $n \to \infty$.  In particular, if $q \to 0$, then it has risk tending to 0.
\end{thm}

In essence, LOND (with a proper choice of $\nu$ above) has risk tending to 0 when $r - (1 - r^{1/\gamma})^\gamma > \beta$.
This is the best upper bound that we were able to establish for the LOND algorithm.  We do not know if it is optimal or not.  In particular, it's quite possible that LOND also achieves the static selection boundary.  

\begin{rem}
The analog of \remref{lambda} applies here as well.  (Technical details are omitted.)
\end{rem}

\section{Numerical experiments}
\label{sec:numerics}

In this section, we perform some simulations to study the performance of LORD and LOND algorithms on finite data, and also to compare them with the (static) BH procedure. We consider the normal model and the double-exponential model.  It is worth repeating that the BH procedure, which is a static procedure, requires knowledge of all P-values to determine the significance level for testing the hypotheses.  Hence, it does not address the scenario in online testing.  In contrast, the sequential methods decide the significance level at each step based on previous outcomes and are required to control de FDR at each step.

In our experiments, for both LORD and LOND, we choose the sequence $(\lambda_i)_{i = 1}^{\infty}$ as
\beq \label{level}
\lambda_i = L i^{-1.05},
\eeq
with $L$ set to ensure $\sum_{i = 1}^{\infty} \lambda_{i} = q$, where (as before) $q$ denotes the desired FDR level.
 
\subsection{Fixed sample size} \label{sec:fixed}
In this first set of experiments, the sample size is chosen large at $n = 10^9$. 
We draw $m$ observations from the alternative distribution $\Phi(\cdot - \mu)$, and the other $n-m$ from the null distribution $\Phi$.  All the models are parameterized as in \eqref{n-mu}. 
We choose a few values for the parameter $\beta$ so as to exhibit different sparsity levels, while the parameter $r$ takes values in a grid of spanning $[0,1.5]$. Each situation is repeated 300 times and we report the average FDP and FNP for each procedure.
The FDR control level is set at $q = 0.1$.  

\subsubsection{Normal model}
\begin{figure}[h!]\centering
	\includegraphics[width = 5cm, height= 4.9cm]{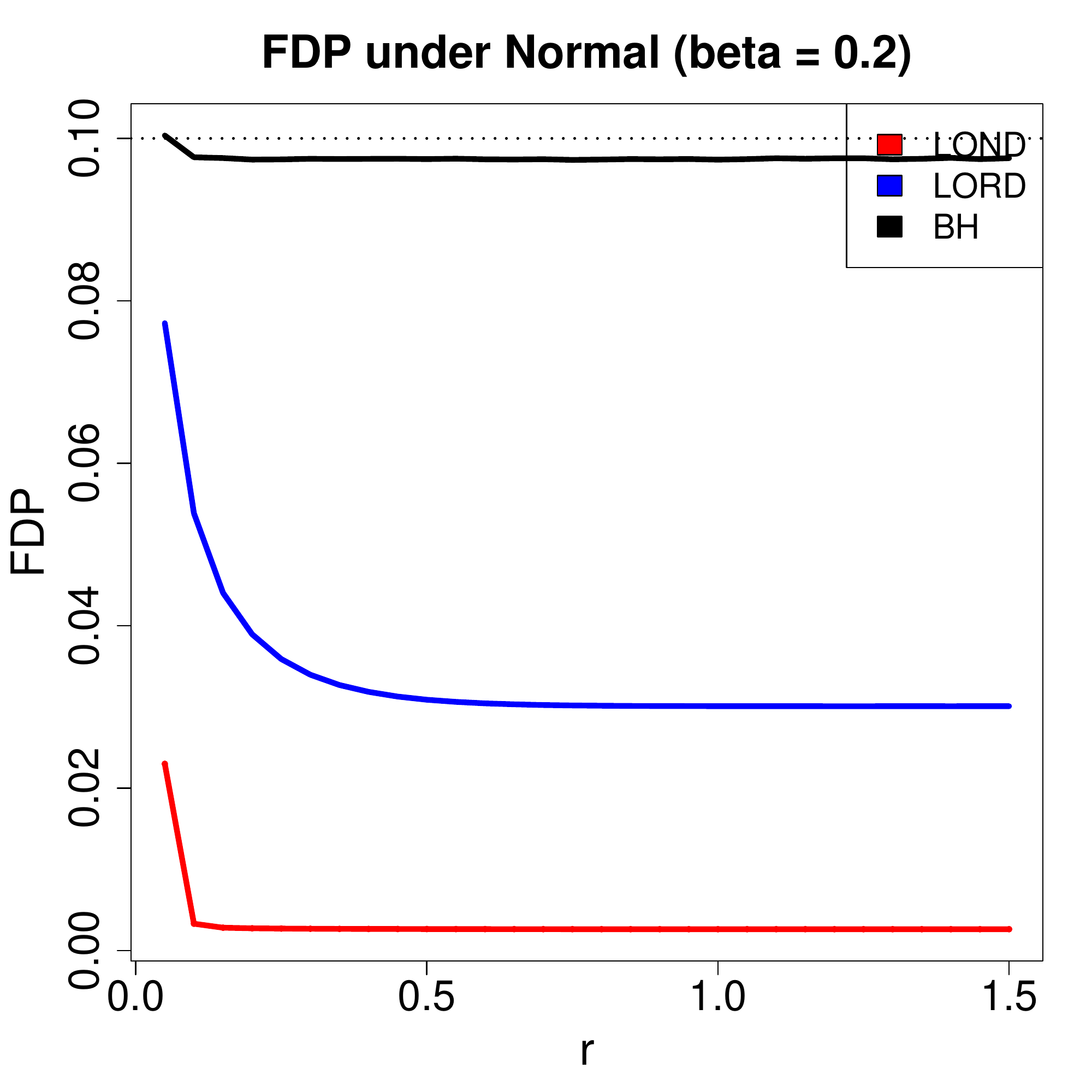}
	\includegraphics[width = 5cm, height= 4.9cm]{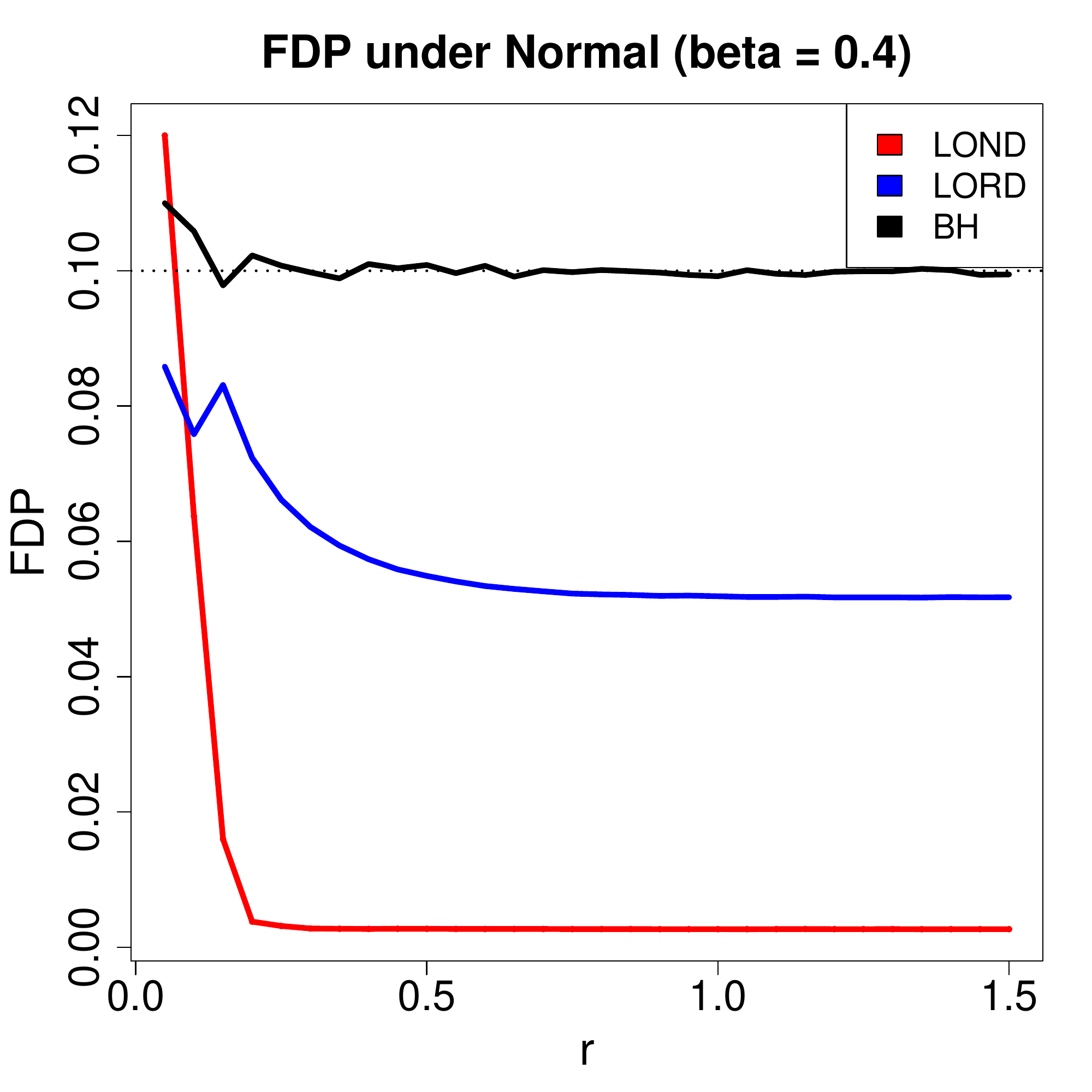}
	\includegraphics[width = 5cm, height= 4.9cm]{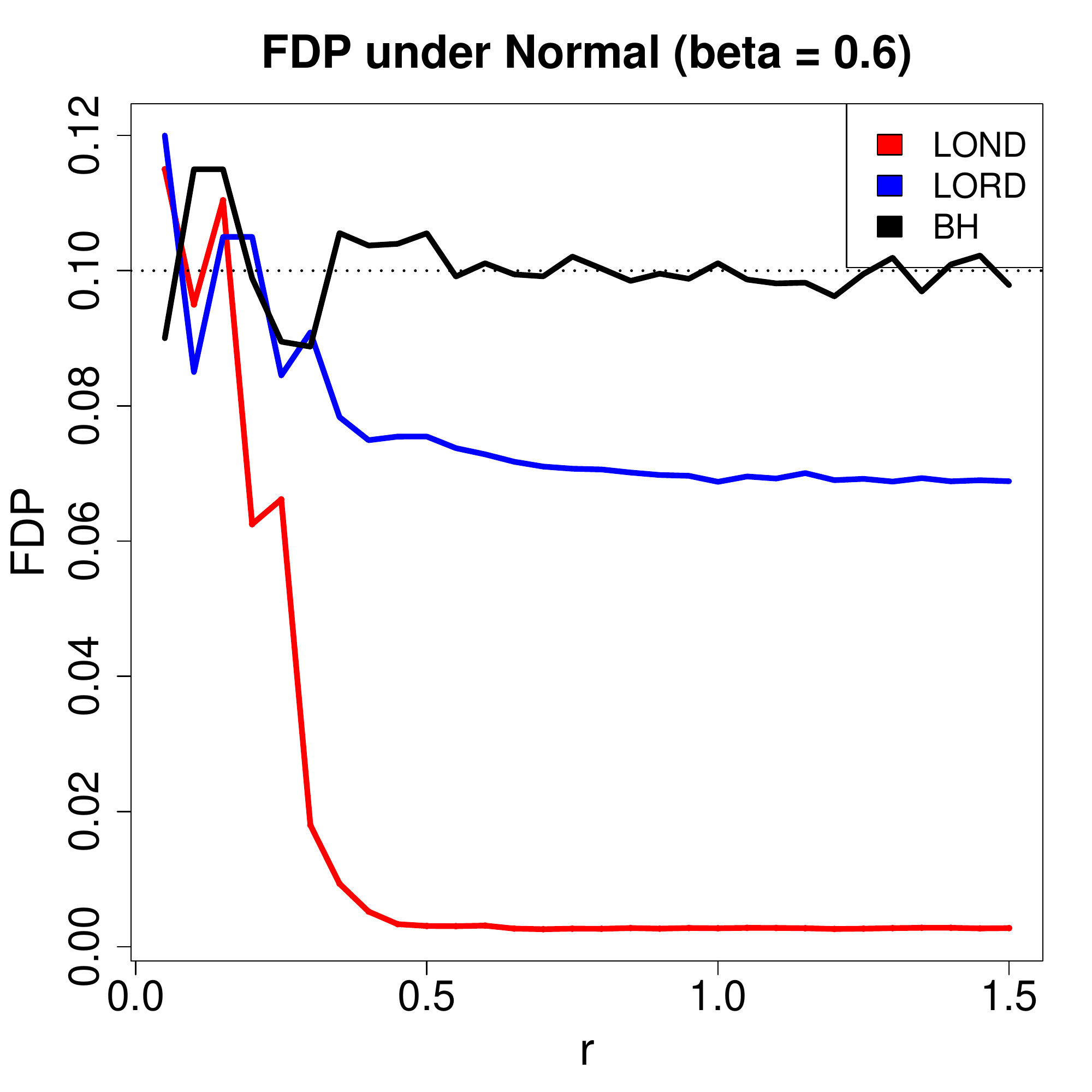}
	\caption{Simulation results showing the FDP for the BH, LORD and LOND methods under the normal model in three distinct sparsity regimes. The black horizontal line delineates the desired FDR control level ($q = 0.1$).}
	\label{fig:fdp_normal}	
\end{figure}

\begin{figure}[h!]\centering
	\includegraphics[width = 5cm, height= 5cm]{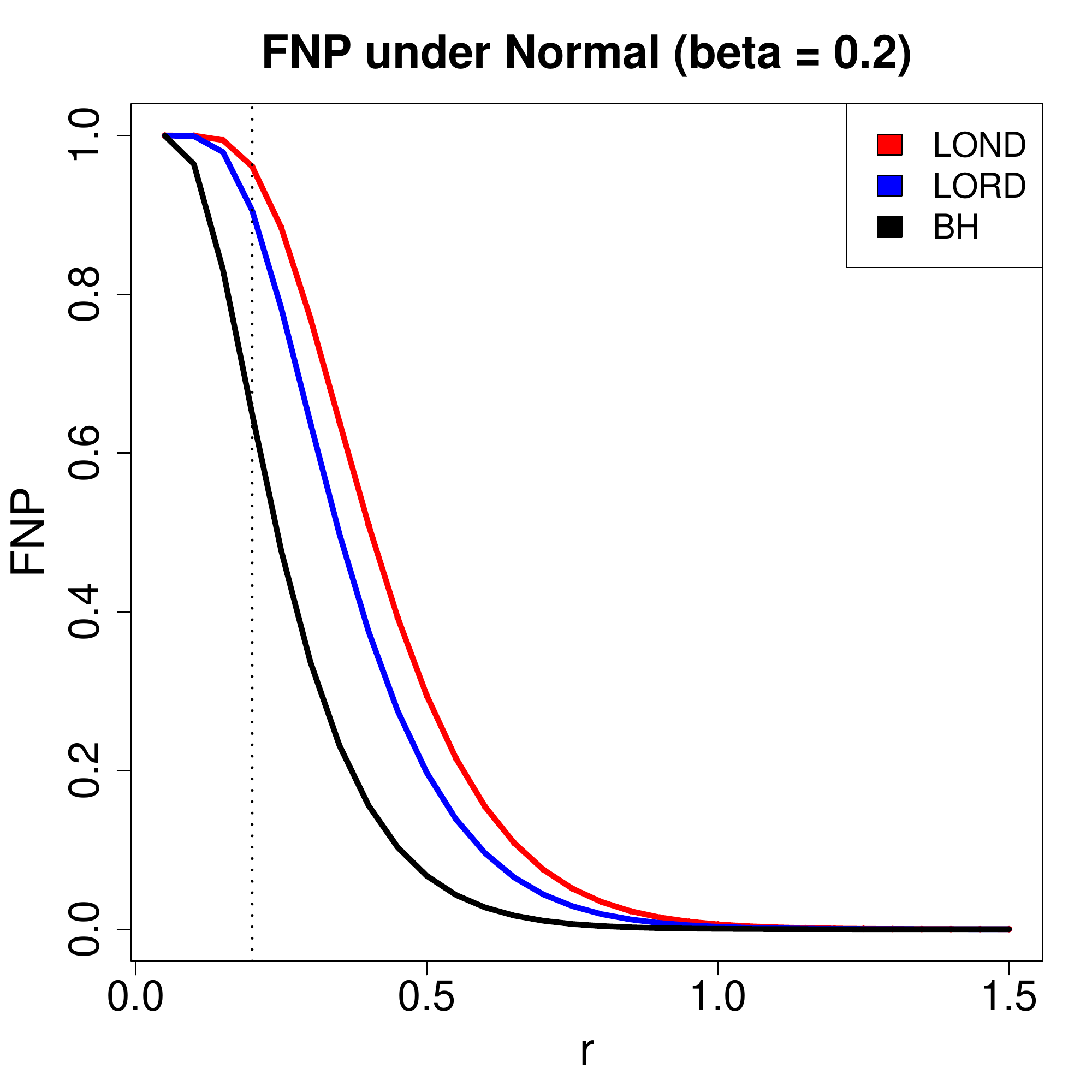}
	\includegraphics[width = 5cm, height= 5cm]{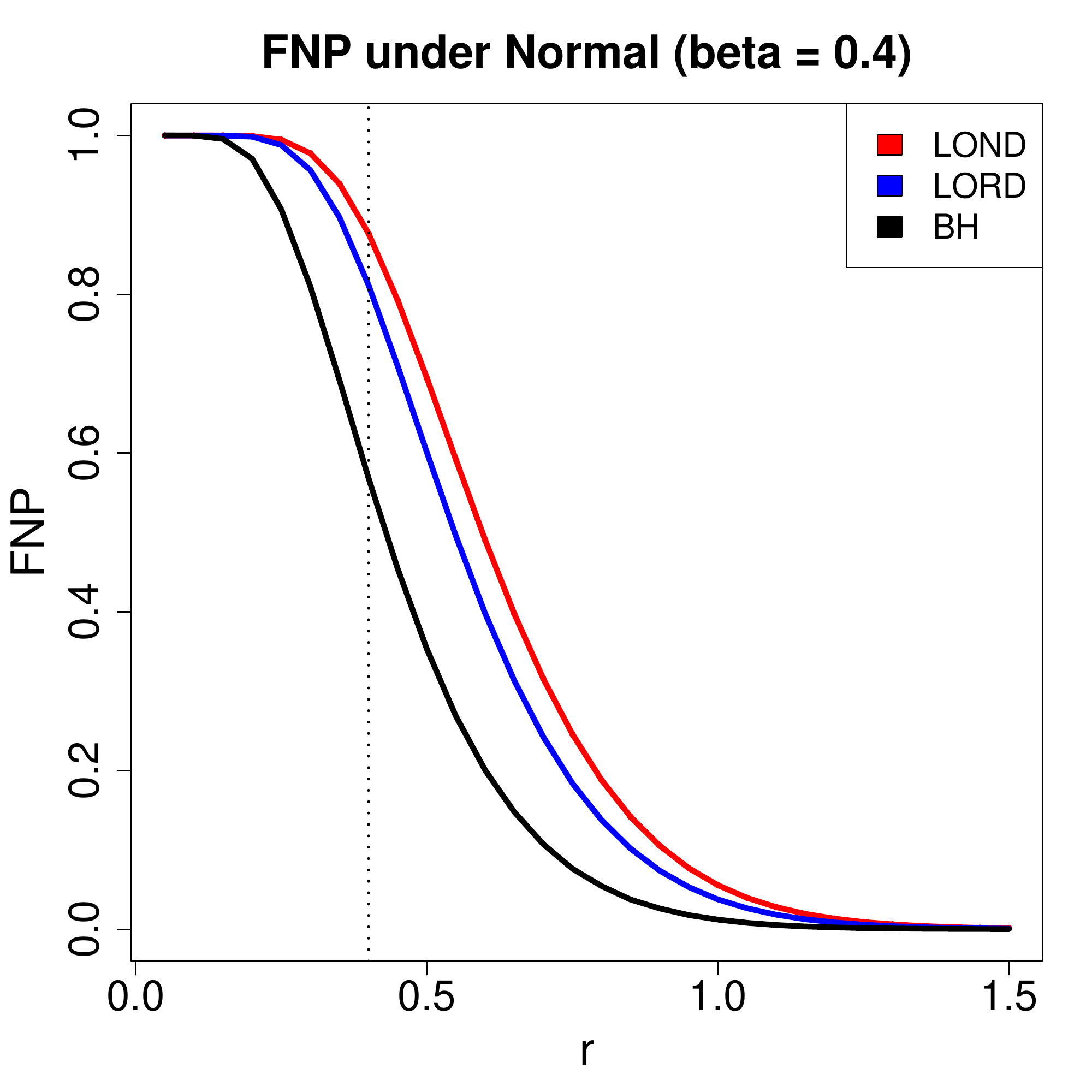}
	\includegraphics[width = 5cm, height= 5cm]{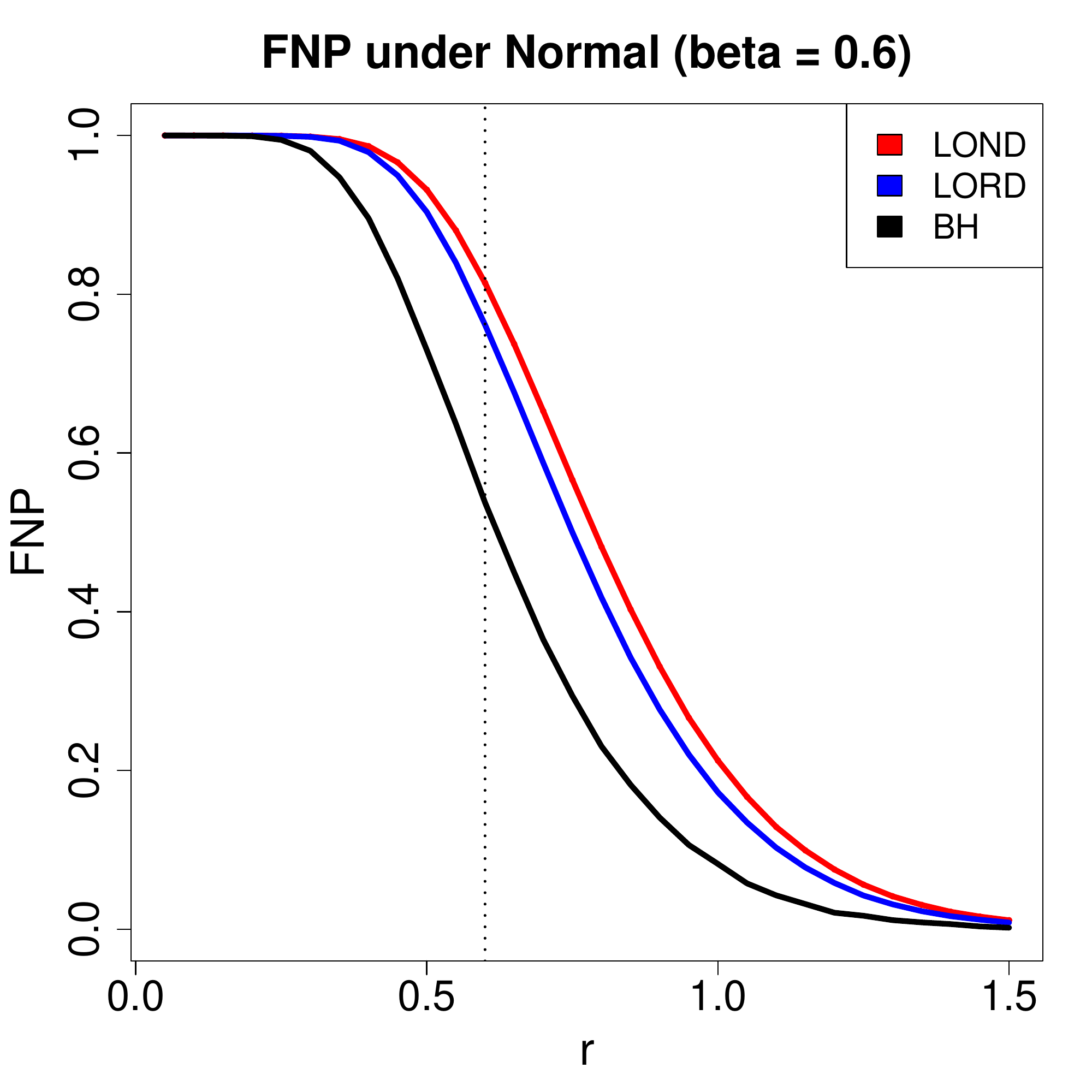}
	\caption{Simulation results showing the FNP for the BH, LORD and LOND methods under the normal model in three distinct sparsity regimes. The black vertical line delineates the theoretical threshold  ($r=\beta$).}
	\label{fig:fnp_normal}
\end{figure}

In this model $\Phi$ is the standard normal distribution. 
The simulation results are reported in \figref{fdp_normal} and \figref{fnp_normal}.
In \figref{fdp_normal} we report the FDP.  We see that LOND becomes more conservative than LORD as $r$ increases.  
In \figref{fnp_normal} we report the FNP.  
We see that LOND is clearly less powerful than LORD in the regime $\beta = 0.2$, but performs comparably to LORD in the regime $\beta = 0.6$ . This is in line with the theory that LOND can at least achieve the line $r = \beta + (1 - r^{1/\gamma})^\gamma$, which is getting closer to $r = \beta$ with increasing values of $\beta$.
We notice that both LORD and LOND are clearly less powerful than BH in finite samples, even at $n = 10^9$, even though our theory says that LORD achieves the same selection boundary as BH in the large-sample limit. Also, due to the limitation in choice of $\nu$ (here $\nu = 1.05$), the selection boundary that LORD can achieve is $r = \nu \beta $ by theory, which explains why LORD lags behind BH.
Finally, we remark the transition of LORD from high FNP to low FNP happens in the vicinity of the theoretical threshold ($r = \beta$).

\subsubsection{Double-exponential model}
\begin{figure}[h!]\centering
	\includegraphics[width = 5cm, height= 4.9cm]{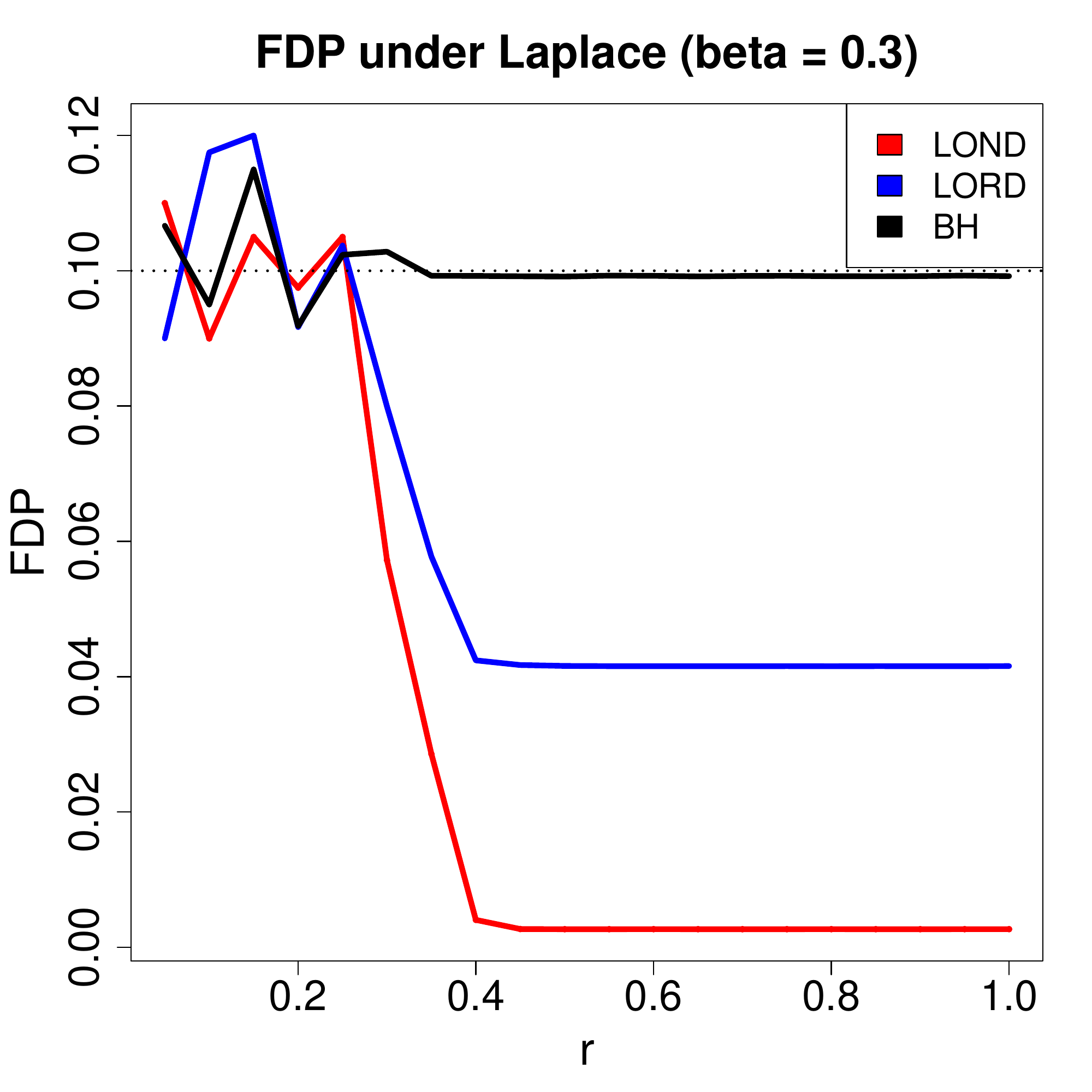}
	\includegraphics[width = 5cm, height= 4.9cm]{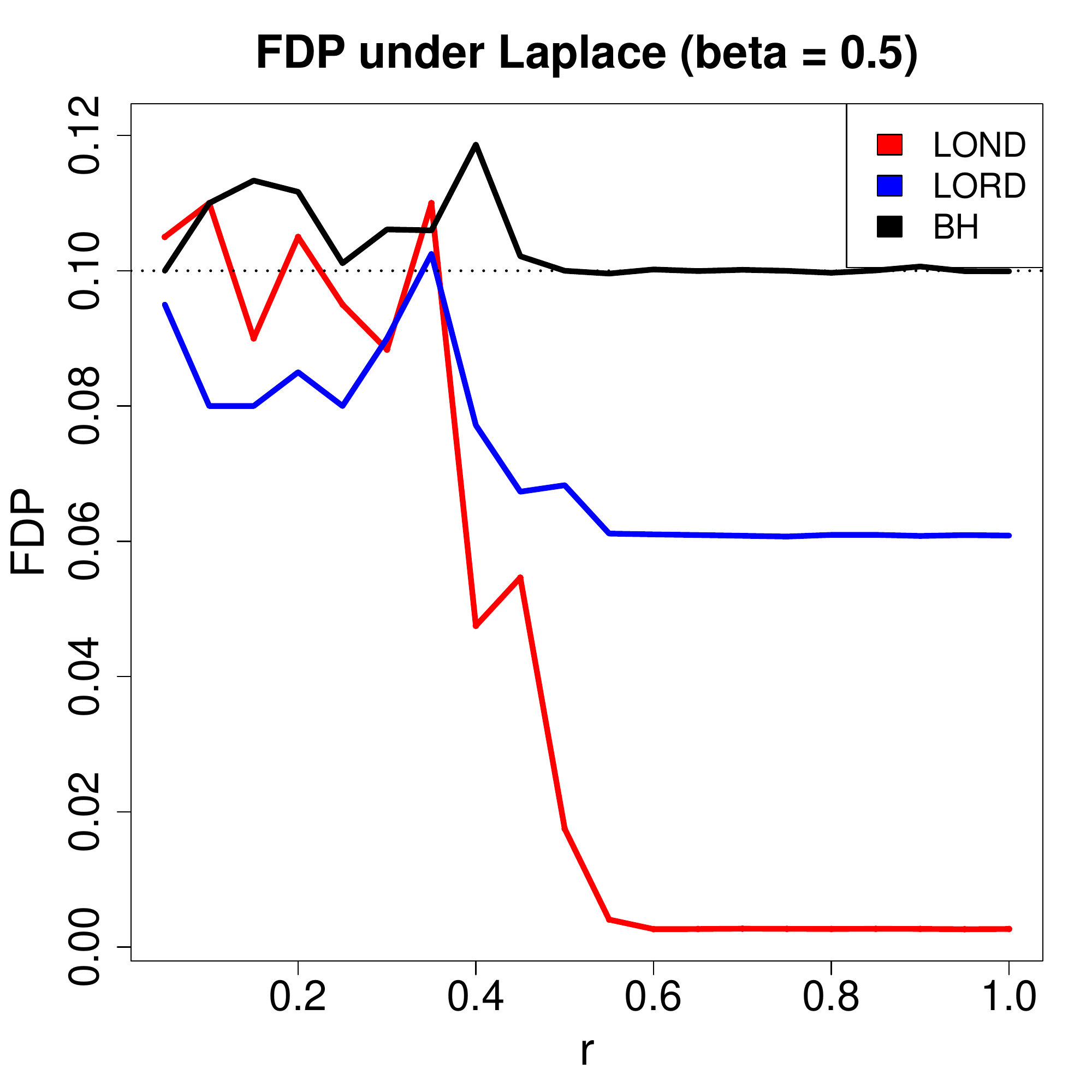}
	\includegraphics[width = 5cm, height= 4.9cm]{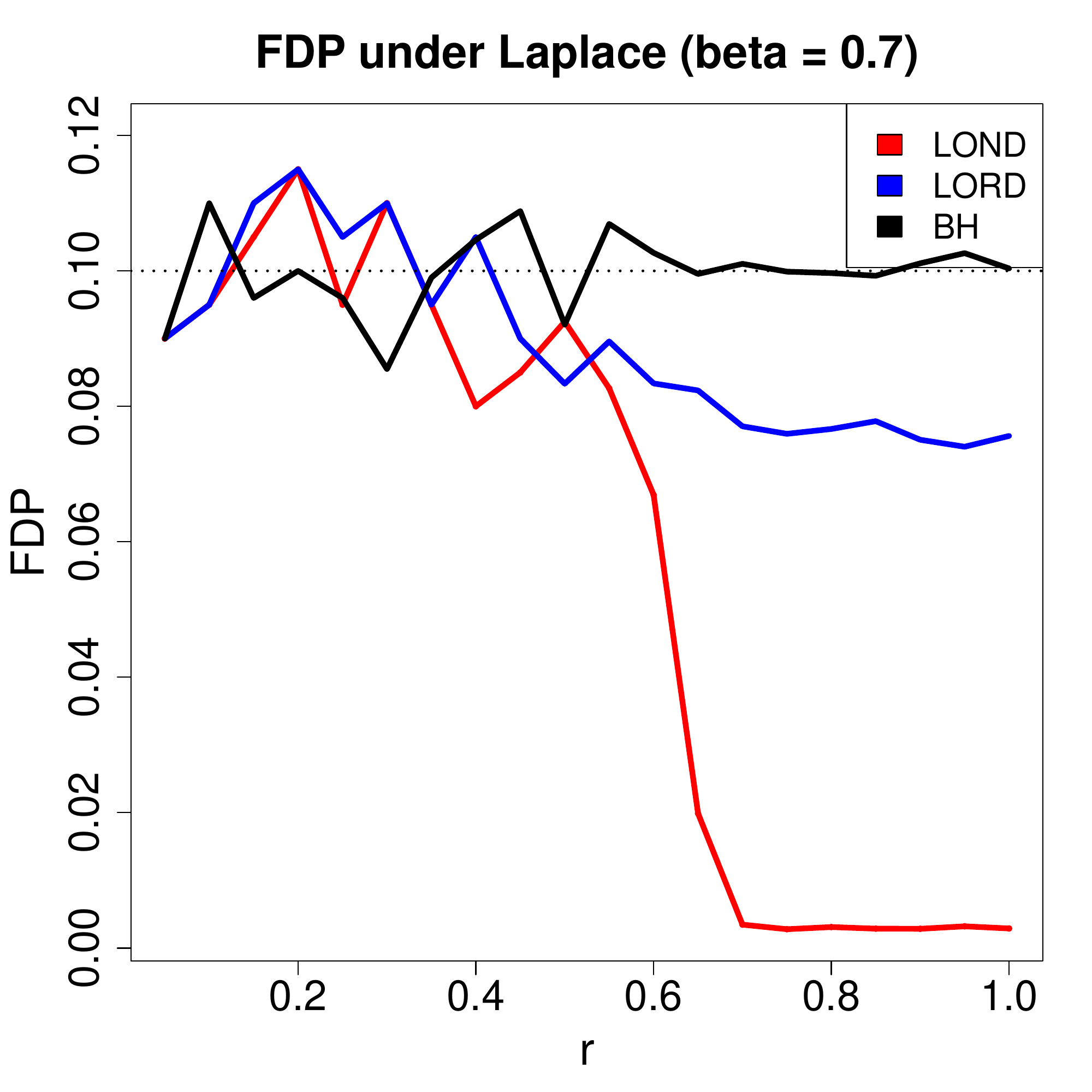}
	\caption{Simulation results showing the FDP for the BH, LORD and LOND methods under the double-exponential model in three distinct sparsity regimes. The black horizontal line delineates the desired FDR control level ($q = 0.1$).}
	\label{fig:fdp_laplace}	
\end{figure}

\begin{figure}[h!]\centering
	\includegraphics[width = 5cm, height= 5cm]{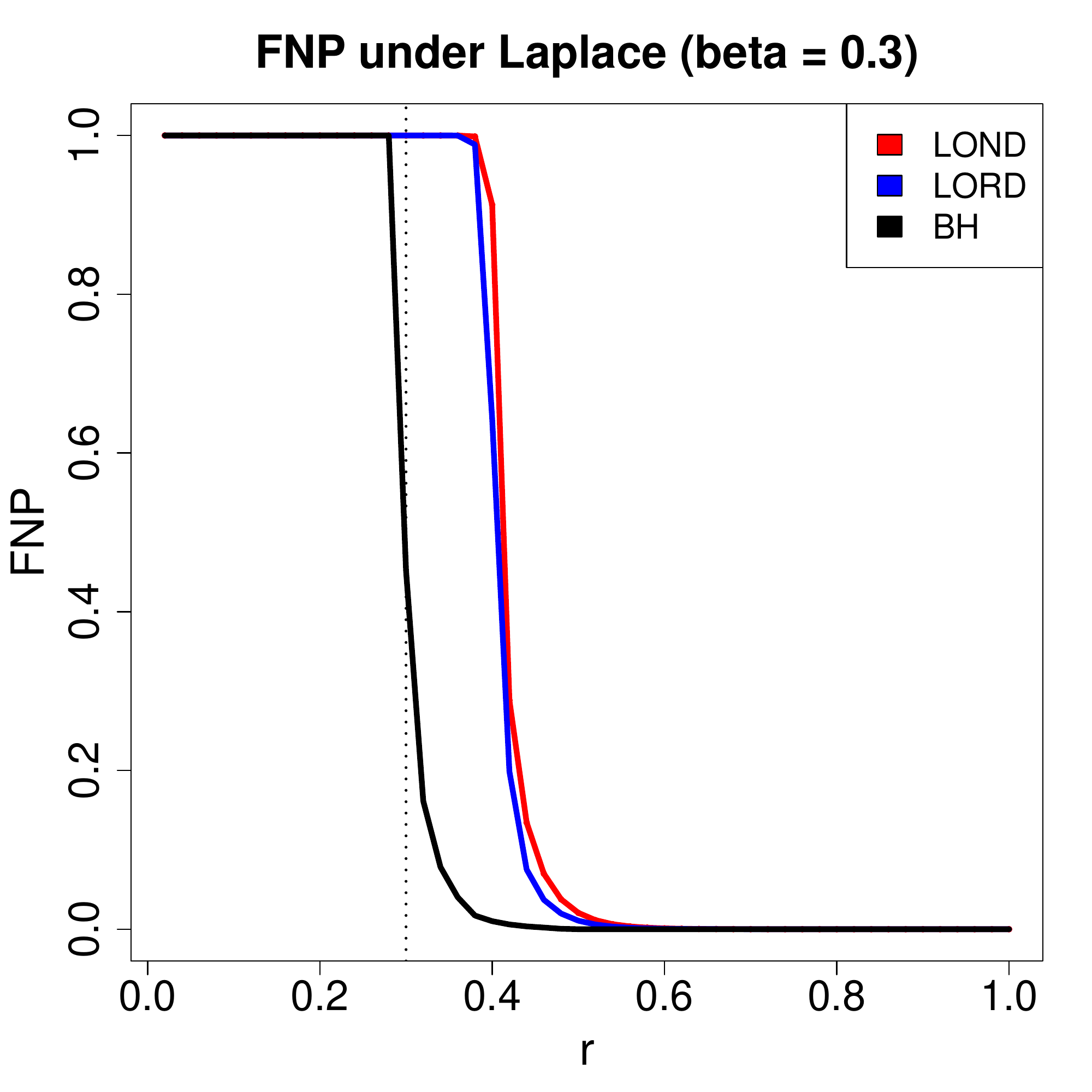}
	\includegraphics[width = 5cm, height= 5cm]{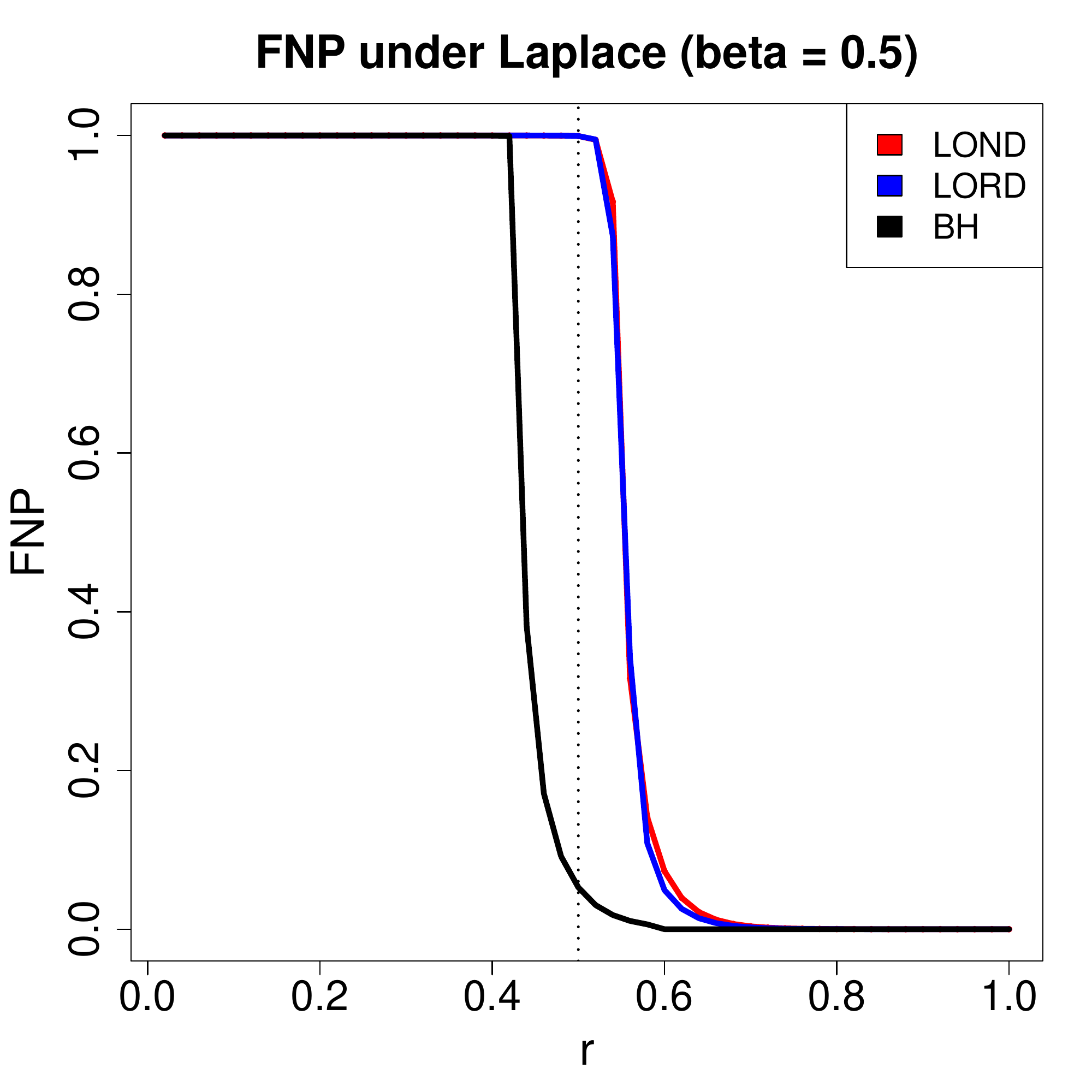}
	\includegraphics[width = 5cm, height= 5cm]{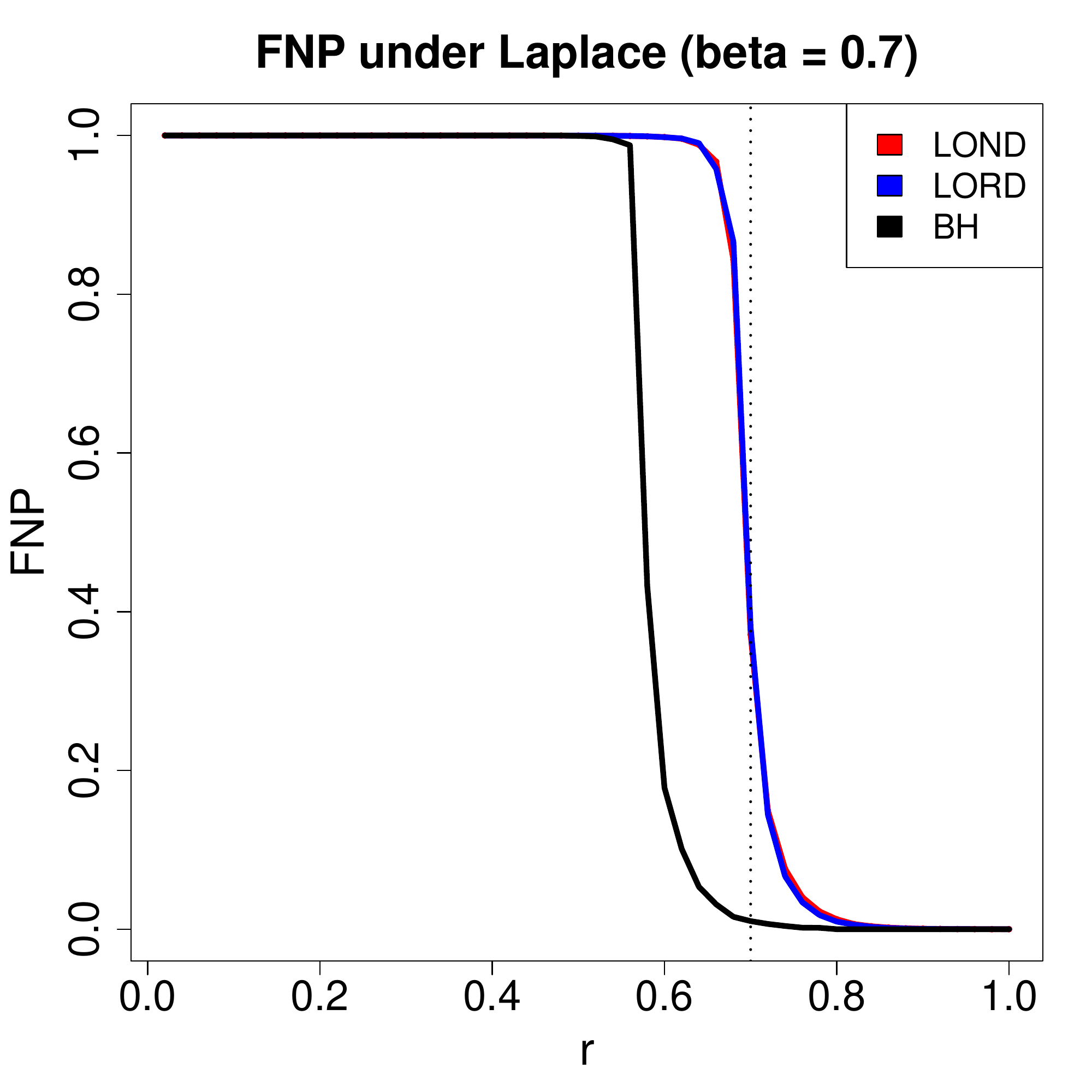}
	\caption{Simulation results showing the FNP for the BH, LORD and LOND methods under the double exponential model in three distinct sparsity regimes. The black vertical line delineates the theoretical threshold  ($r=\beta$).}
	\label{fig:fnp_laplace}
\end{figure}

In this model $\Phi$ is the double-exponential distribution with variance 1. The simulation results are reported in \figref{fdp_laplace} (FDP) and \figref{fnp_laplace} (FNP).
Here we observe that LOND becomes more conservative than LORD as $r$ increases in terms of FDP. The LOND and LORD perform more  comparably than in the normal setting in terms of FNP, especially when $\beta$ is close to 1. This is again in line with our theoretical results. The BH method clearly outperforms the other two methods even though $n = 10^9$.  Due to the limitation in choice of $\nu$ (here $\nu = 1.05$), the selection boundary that LORD can achieve is $r = \nu \beta $ by theory, which explains why LORD lags behind BH.
The transition of three methods from FNP near 1 to FNP near 0 happens, again, in the vicinity of the theoretical threshold, but is much sharper here. 

\subsection{Varying sample size} \label{sec:fixed}
In this second set of experiments we examine the effect of various sample sizes on the risk of the LORD and LOND procedures under the standard normal model and the double-exponential model (with variance~1). 

\subsubsection{FNR of LORD with a fixed level}
\begin{figure}[h!]\centering
	\includegraphics[width = 5cm, height= 5cm]{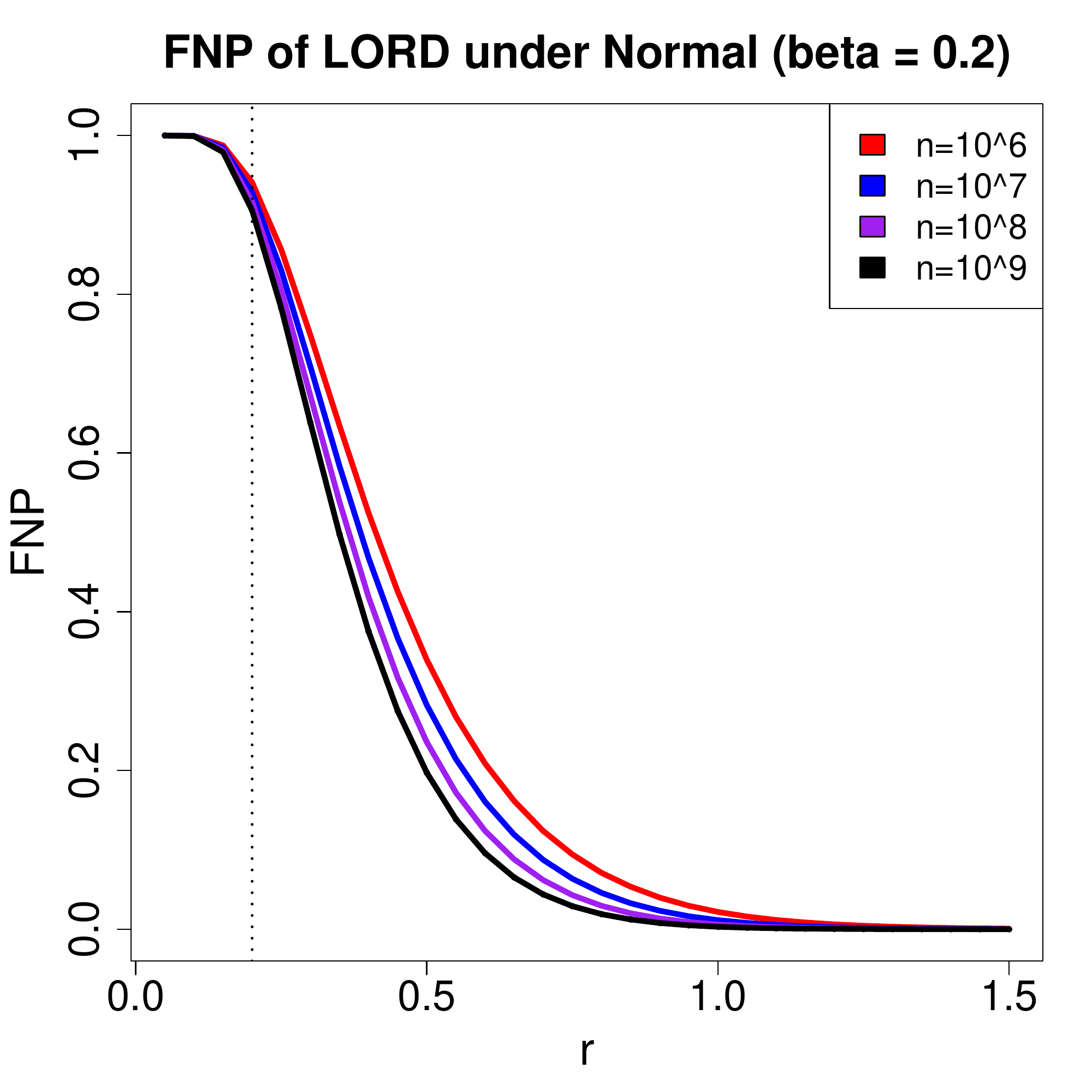}
	\includegraphics[width = 5cm, height= 5cm]{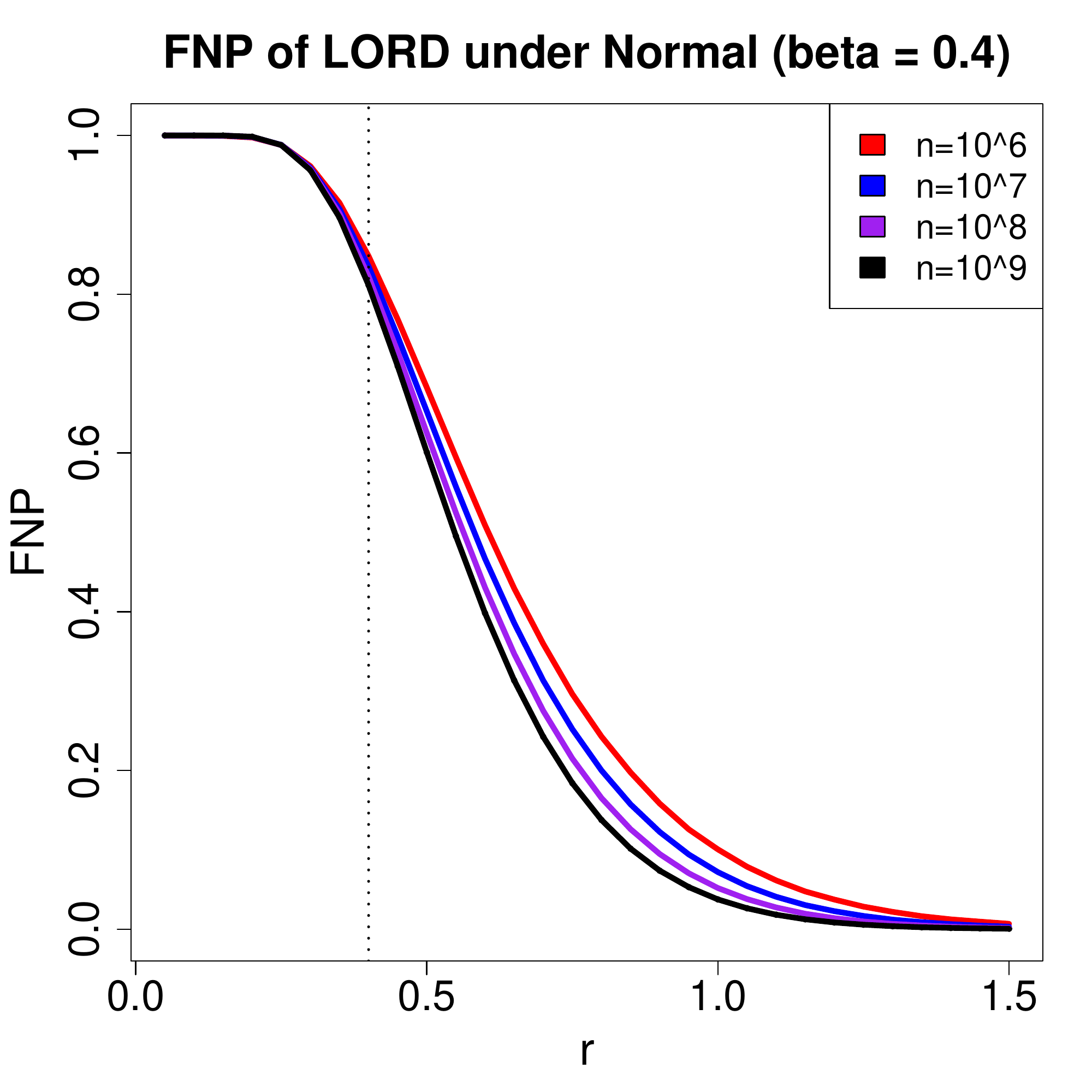}
	\includegraphics[width = 5cm, height= 5cm]{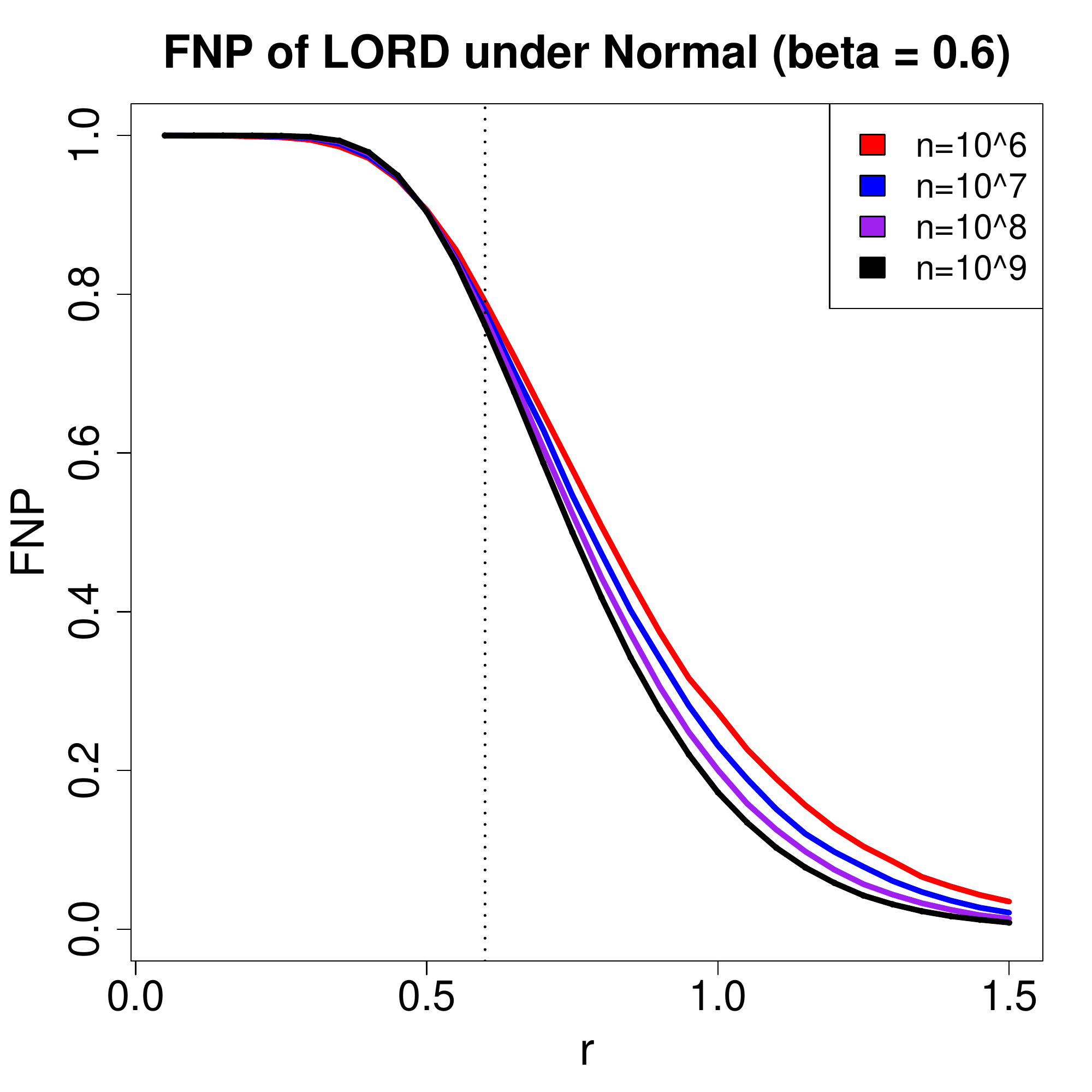}
	\caption{Simulation results showing the FNP for LORD under the normal model in three distinct sparsity regimes with different sample size. The black vertical line delineates the theoretical threshold  ($r=\beta$).}
	\label{fig:fnp_lord_normal}
\end{figure}

\begin{figure}[h!]\centering
	\includegraphics[width = 5cm, height= 5cm]{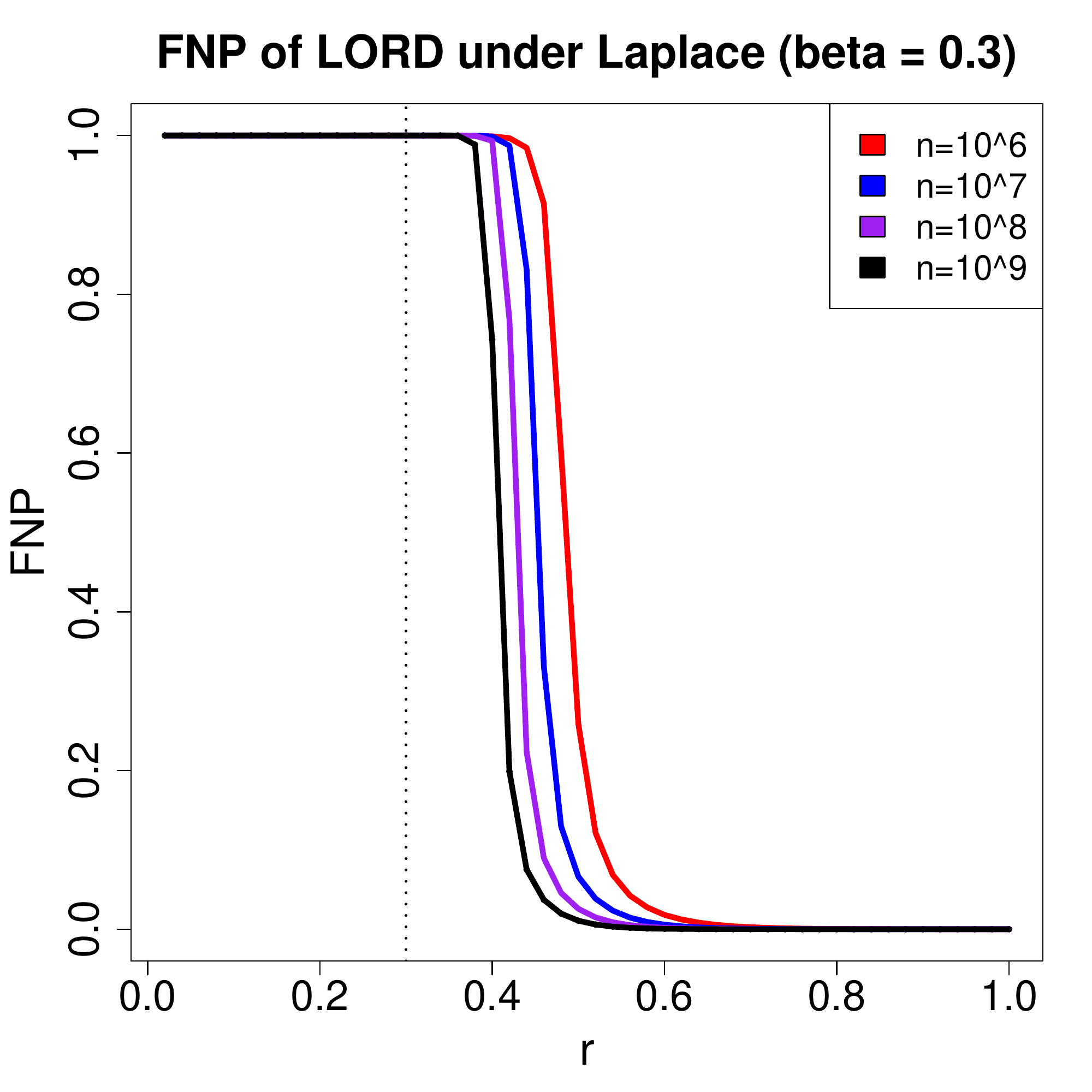}
	\includegraphics[width = 5cm, height= 5cm]{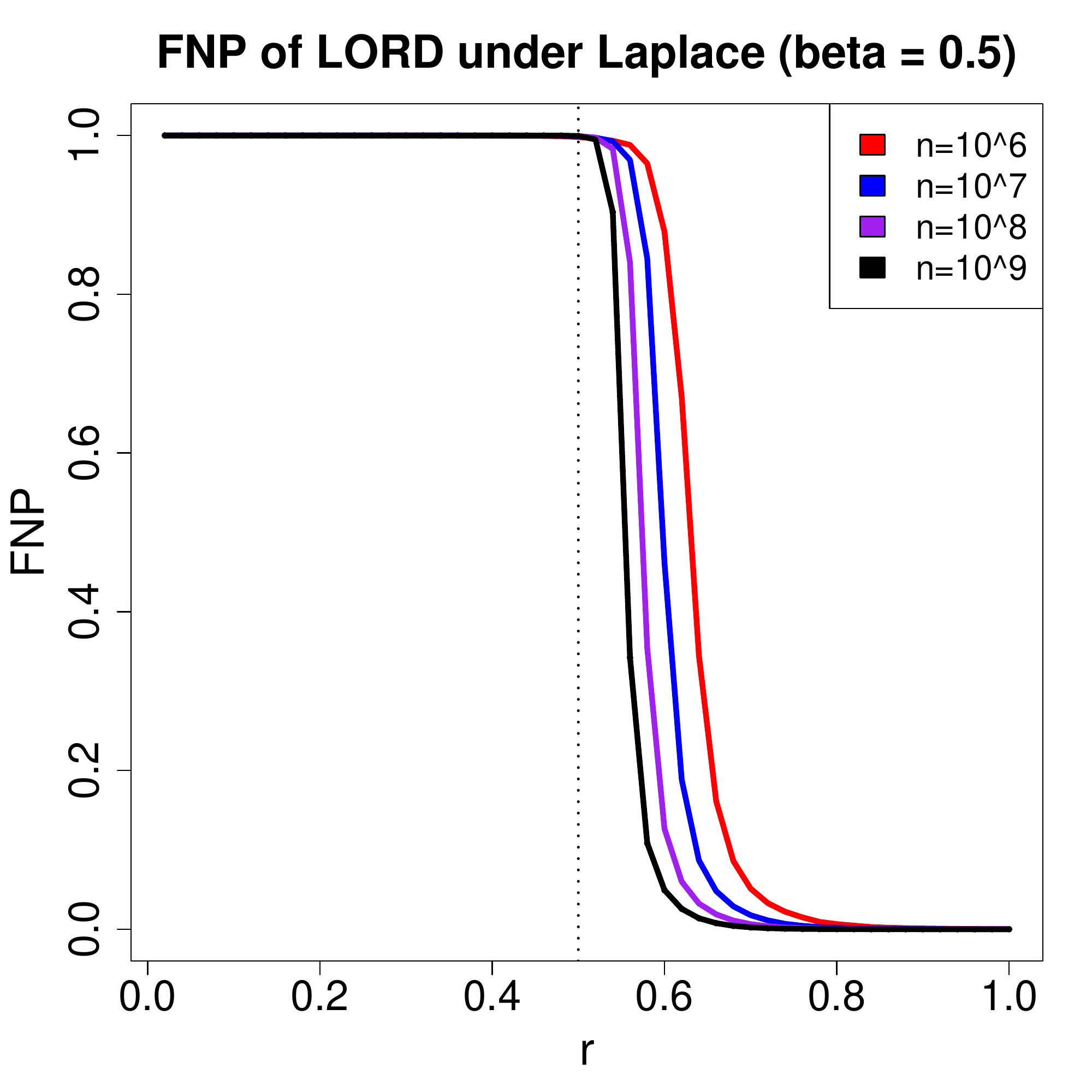}
	\includegraphics[width = 5cm, height= 5cm]{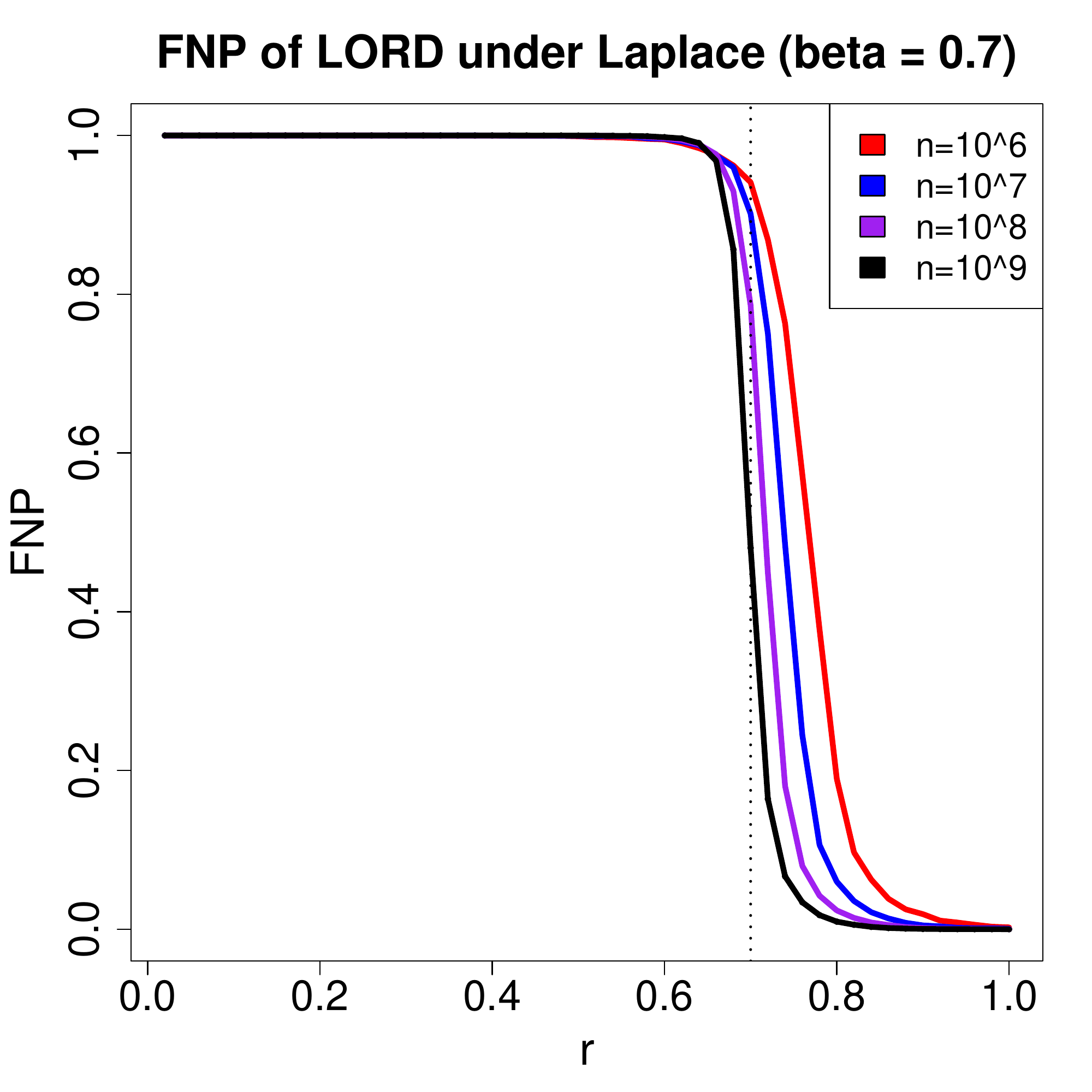}
	\caption{Simulation results showing the FNP for LORD under the double exponential model in three distinct sparsity regimes with different sample size. The black vertical line delineates the theoretical threshold  ($r=\beta$).}
	\label{fig:fnp_lord_laplace}
\end{figure}

In this subsection, we present numerical experiments meant to illustrate the theoretical results we derived about asymptotic FNR of LORD. We fix $q = 0.1$ , and choose a few values for the parameter $\beta$ so as to exhibit different sparsity levels, while the parameter $r$ takes values in a grid of spanning $[0,1.5]$.  
We plot the average FNP of LORD procedure with different $n \in \{10^6, 10^7, 10^8, 10^9\}$. The simulation results are reported in \figref{fnp_lord_normal} and \figref{fnp_lord_laplace}. Each situation is repeated 200 times. We observe that in the normal model when $r>\beta$, the FNP decreases as n is getting larger. In the double-exponential model, as $n$ increases, the FNP transition lines are getting closer the theoretical thresholds $r = \beta$, especially when $\beta = 0.7$.

\subsubsection{Varying level}
Here we explore the effect of letting the desired FDR control level $q$ tend to 0 as $n$ increases in accordance with \eqref{q}.  Specifically, we set it as $q = q_n = 1/\log n$.  We choose $n$ on a log scale, specifically, $n \in \{10^5, 10^6, 10^7, 10^8, 10^9\}$.  Each time, we fix a value of $(\beta, r)$ such that $r > \beta$. 

\begin{figure}[h!]\centering
	\includegraphics[width = 5cm, height=5cm]{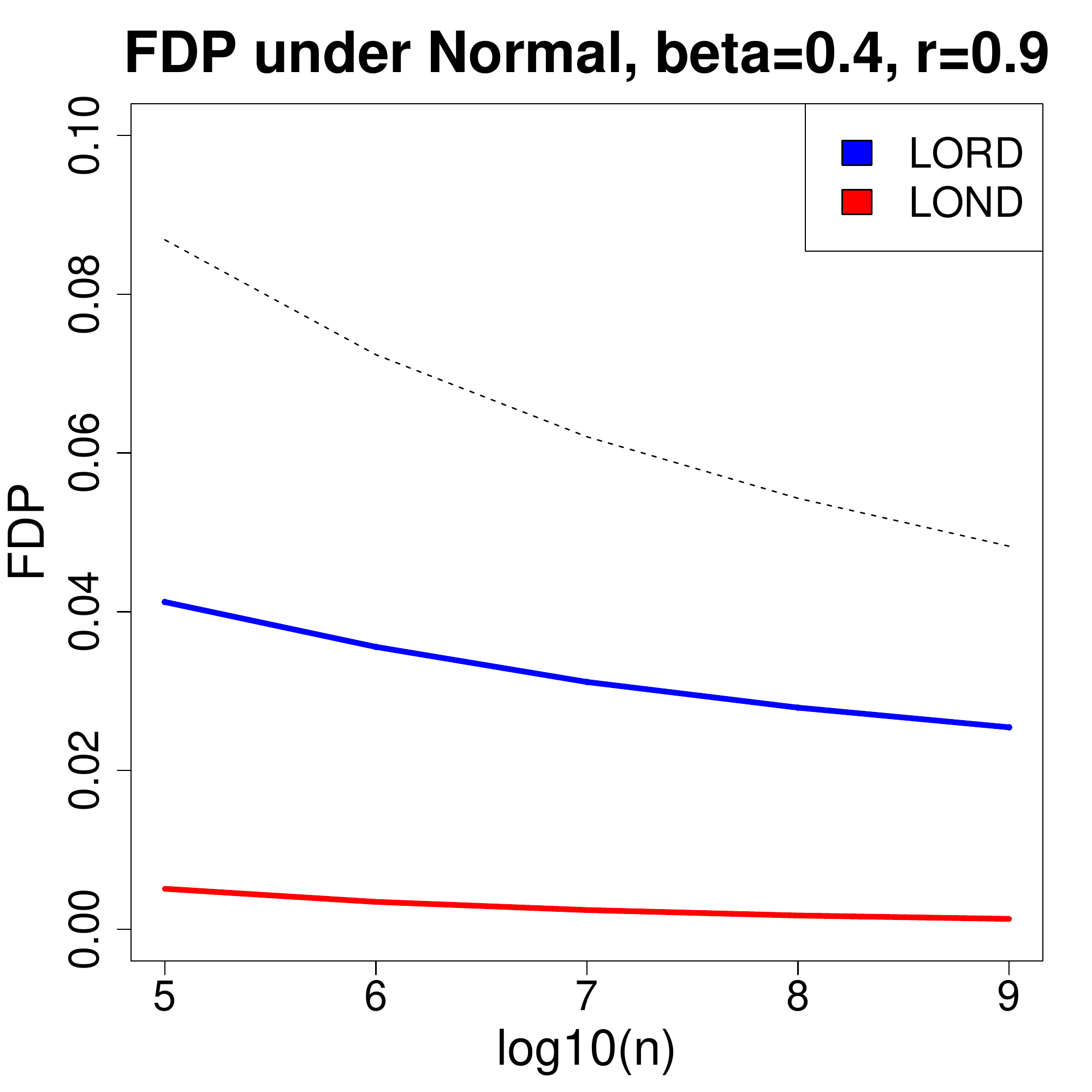}
	\includegraphics[width = 5cm, height=5cm]{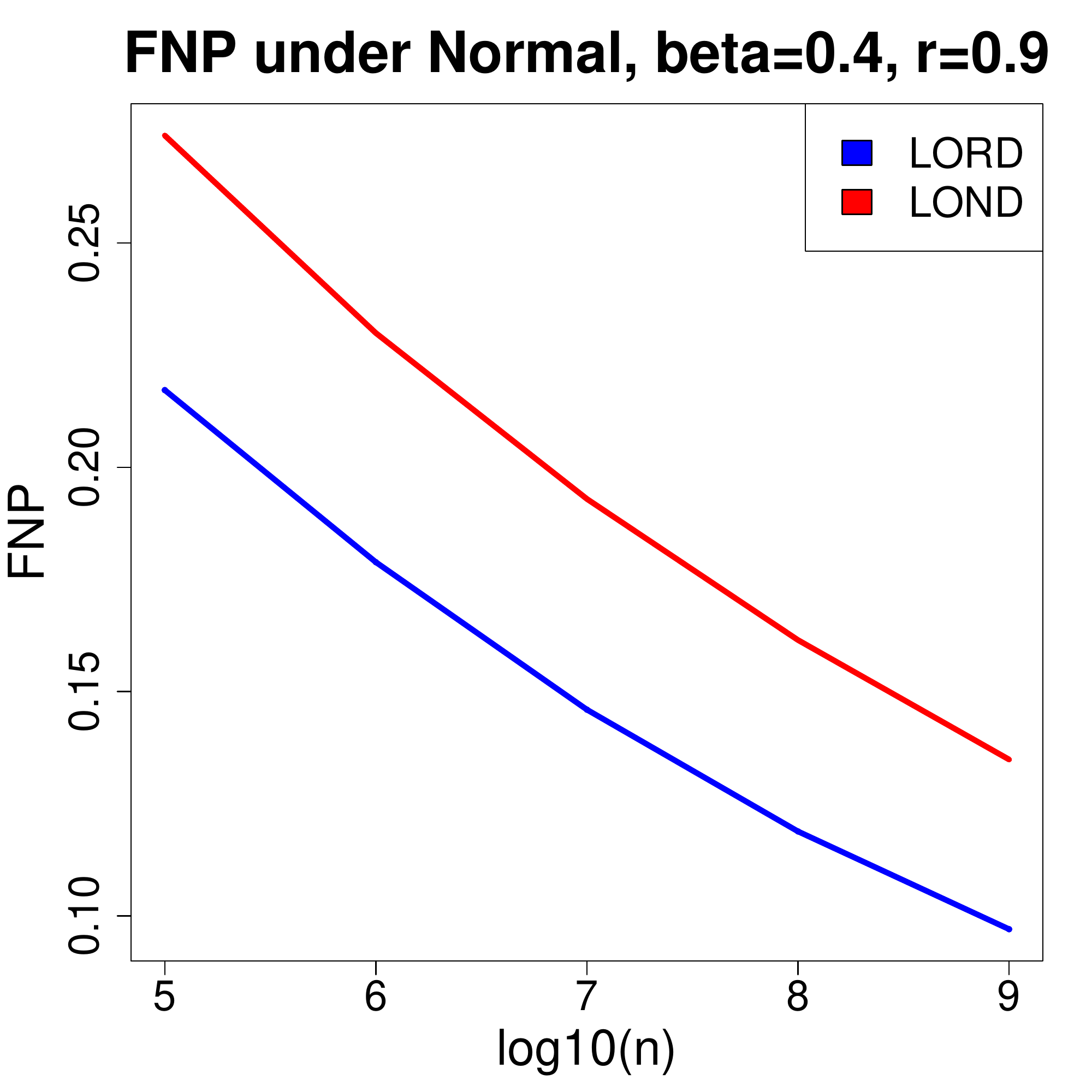}
	\caption{FDP and FNP for the LORD and LOND methods under the normal model with $(\beta, r) = (0.4, 0.9)$ and varying sample size $n$. The black line delineates the desired FDR control level ($q = q_n $).}
	\label{fig:vary_dense_normal}
\end{figure}

\begin{figure}[h!]\centering
	\includegraphics[width = 5cm, height= 5cm]{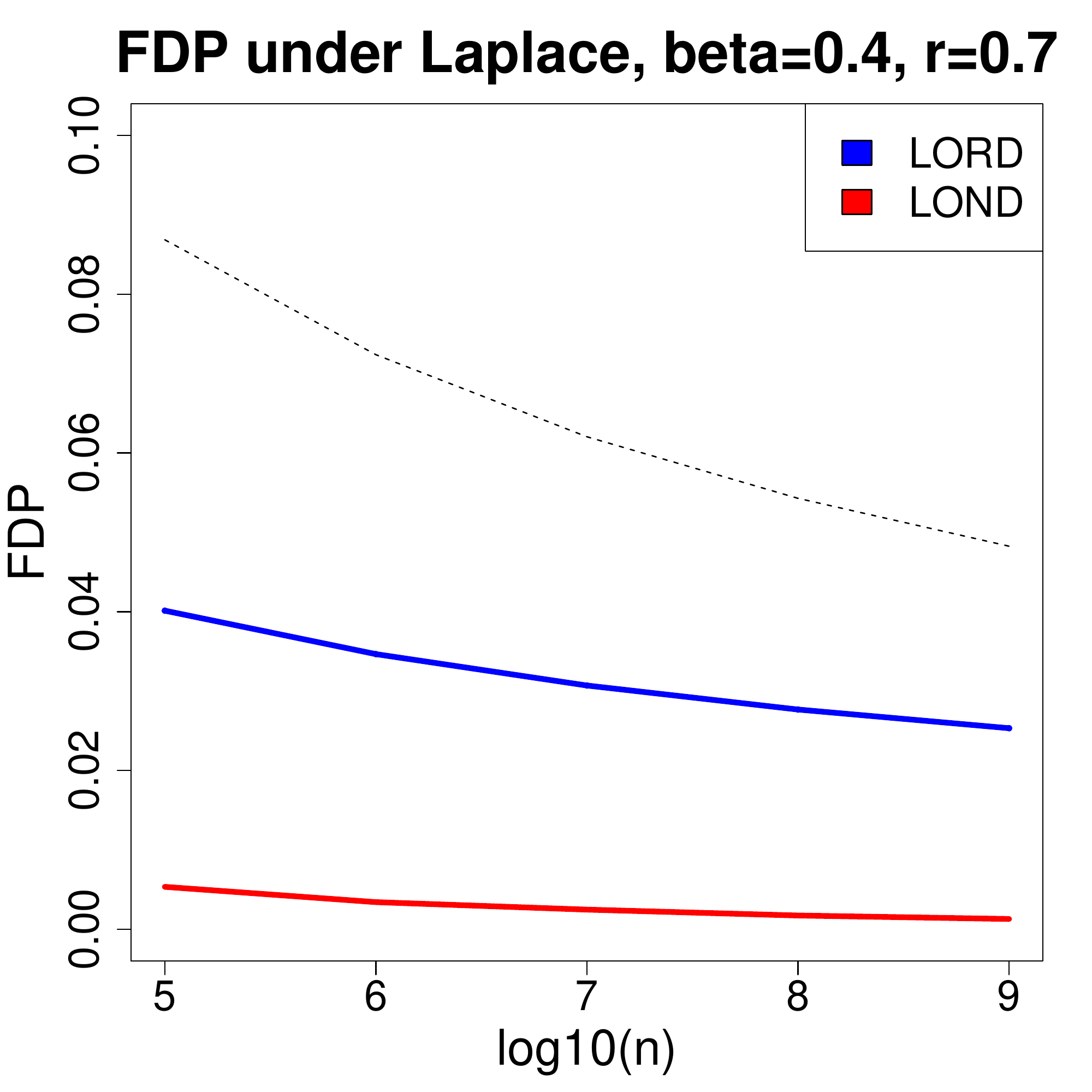}
	\includegraphics[width = 5cm, height= 5cm]{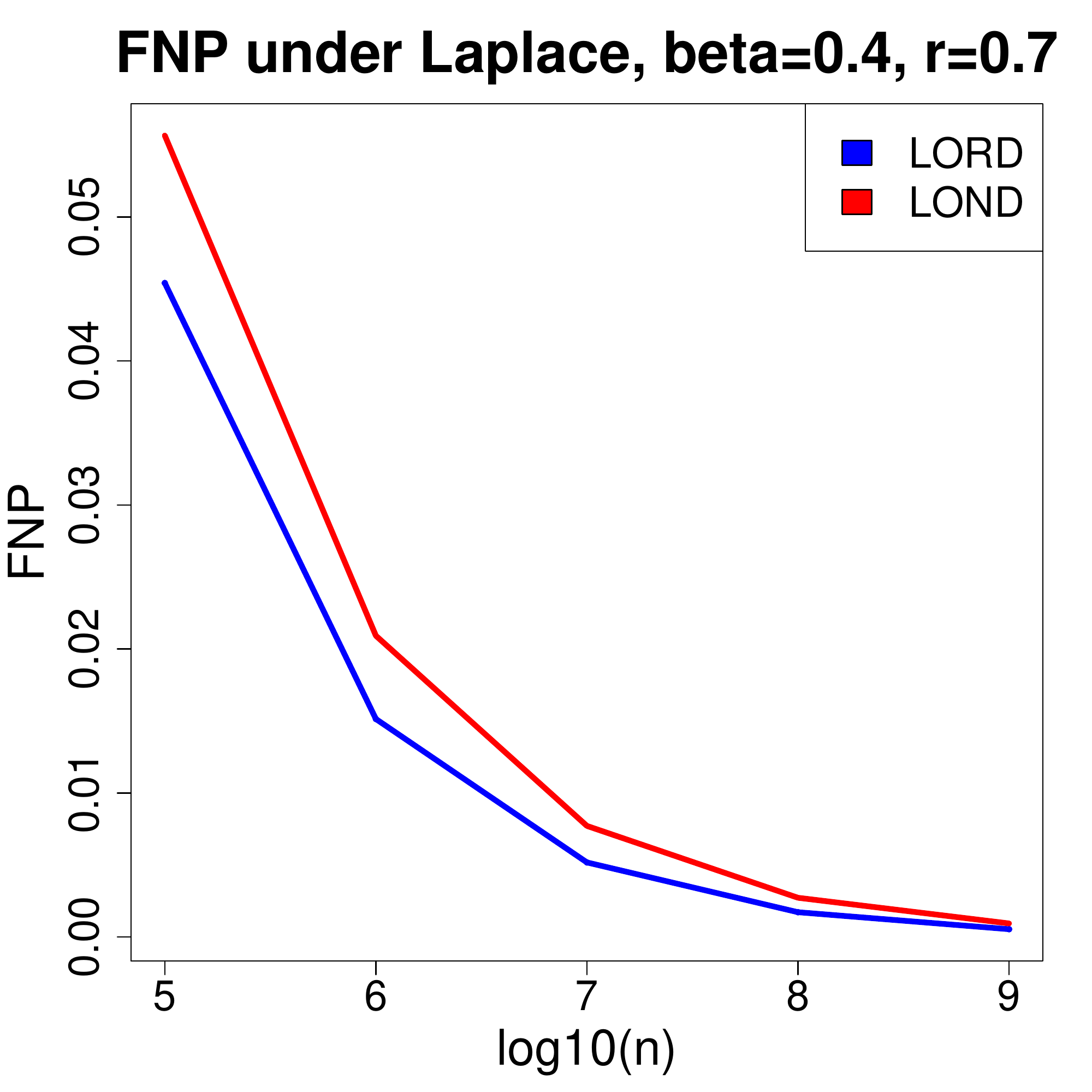}
	\caption{FDP and FNP for the LORD and LOND methods under the double-exponential model with $(\beta, r) = (0.4, 0.7)$ and varying sample size $n$. The black line delineates the desired FDR control level ($q = q_n $). }
	\label{fig:vary_dense_laplace}
\end{figure}

In the first setting, we set $(\beta, r) = (0.4, 0.9)$ for normal model and $(\beta, r) = (0.4, 0.7)$ for double-exponential model.
The simulation results are reported in \figref{vary_dense_normal} and \figref{vary_dense_laplace}.
We see that, in both models, the risks of the two procedures decrease to zero as the sample size gets larger.  LORD clearly dominates LOND (in terms of FNP). 
Both methods have FDP much lower than the level $q_n$, and in particular, LOND is very conservative.

\begin{figure}[h!]\centering
	\includegraphics[width = 5cm, height=5cm]{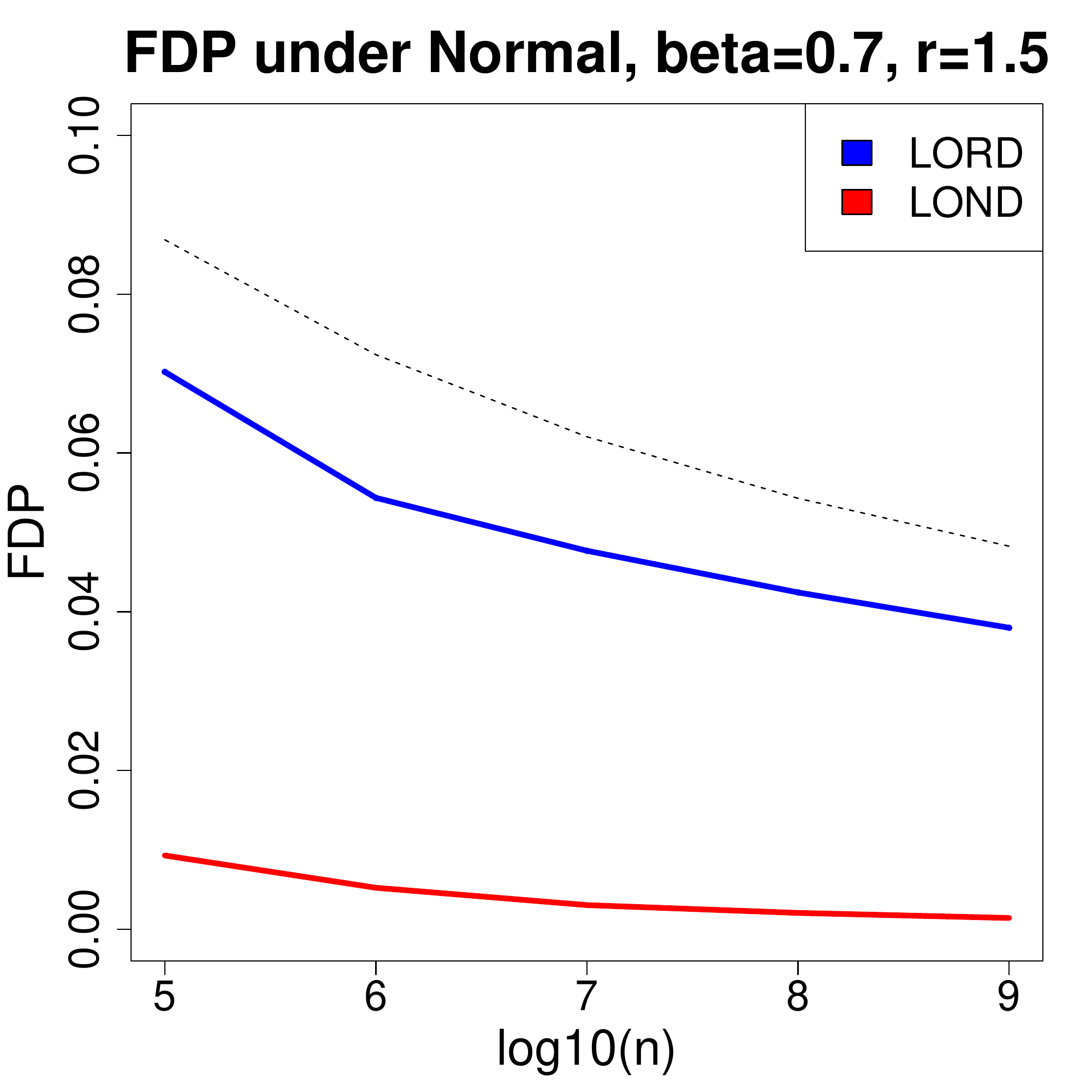}
	\includegraphics[width = 5cm, height=5cm]{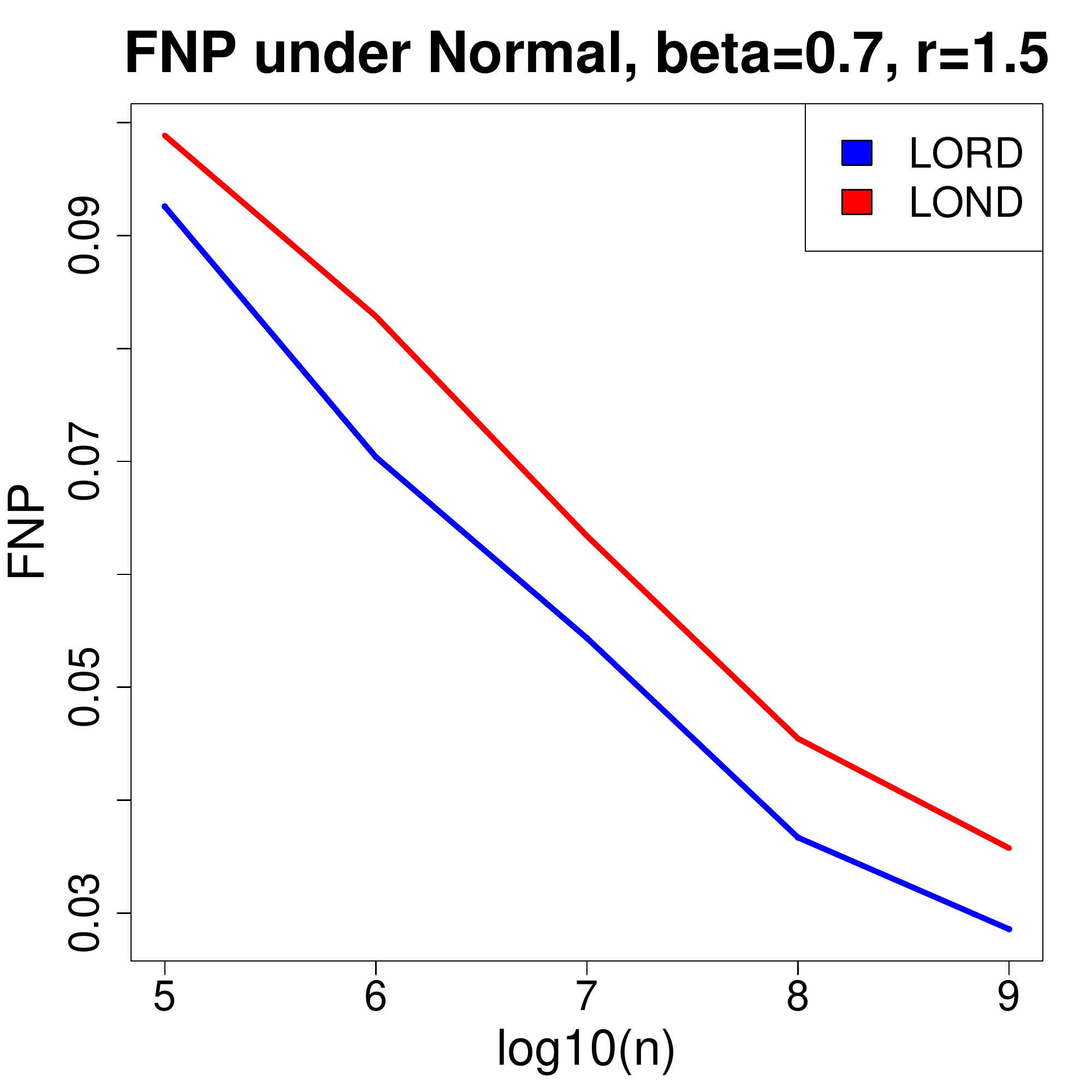}
	\caption{FDP and FNP for the LORD and LOND methods under the normal model with $(\beta, r) = (0.7, 1.5)$ and varying sample size $n$. The black line delineates the desired FDR control level ($q = q_n $).}
	\label{fig:vary_sparse_normal}
\end{figure}

\begin{figure}[h!]\centering
	\includegraphics[width = 5cm, height= 5cm]{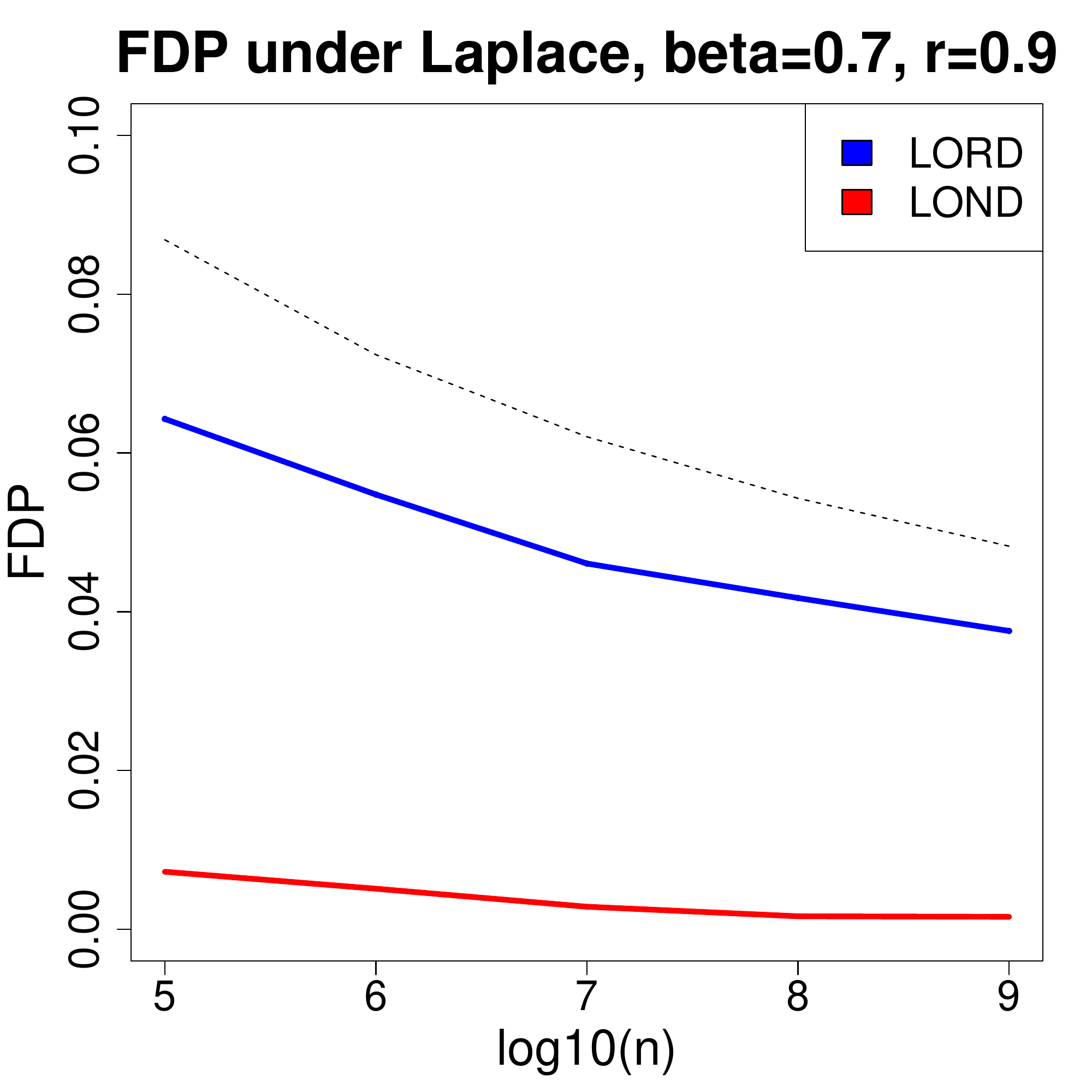}
	\includegraphics[width = 5cm, height= 5cm]{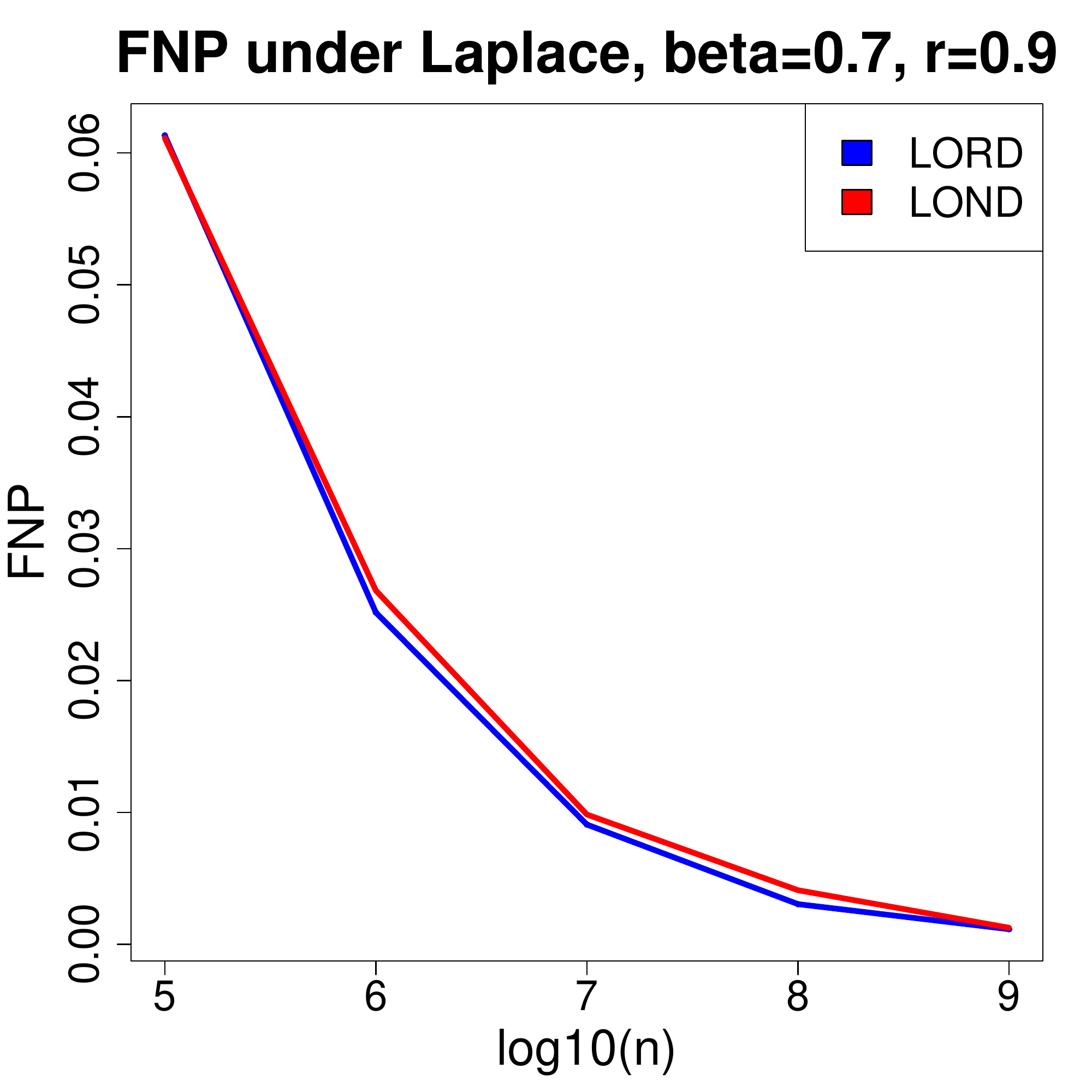}
	\caption{FDP and FNP for the LORD and LOND methods under the double-exponential model with $(\beta, r) = (0.7, 0.9)$ and varying sample size $n$.}
	\label{fig:vary_sparse_laplace}
\end{figure}

In the second setting, we set $(\beta, r) = (0.7, 1.5)$ for normal model and $(\beta, r) = (0.7, 0.9)$ for double-exponential model. 
The simulation results are reported in \figref{vary_sparse_normal} and \figref{vary_sparse_laplace}.
In this sparser regime, we can see that although LORD still dominates, the difference in FNP between two methods is much smaller than that in dense regime, especially in the double-exponential model. 
Both methods have FDP lower than the level $q_n$, and in particular, LOND is very conservative.


\section{Proofs}
\label{sec:proofs}

We prove our results in this section.
Let $\Phi$ denote the CDF of null distribution.  Without loss of generality, we assume throughout that $\Phi(0) = 1/2$. 
Let $F(t)$ denote the CDF of the P-values under alternatives so that 
\beq\label{F}
F(t) = \Phi(\mu - \Phi^{-1}(1-t)),
\eeq
where $\Phi^{-1}$ is the inverse function of $\Phi$.
Let
\beq\label{G}
G(t) = (1-\eps)t + \eps F(t),
\eeq 
which is the CDF of the P-values from the mixture model.
Let $\bar{F} = 1 - F$, which is the survival function of the P-values under alternatives.  Note that 
\beq
\bar{F}(t) = 1-F(t) = 1 - \Phi(\mu - \xi) = \bar{\Phi}(\mu - \xi), 
\eeq
where $\xi := \Phi^{-1}(1-t)$, or equivalently, $t = \bar{\Phi}(\xi)$. 
Because $\Phi$ is as in \defref{AGG}, when $\xi \to \infty$, we have 
\beq
t = \bar{\Phi}(\xi) = \exp\Big\{-\frac{\xi^\gamma}{\gamma} (1+o(1))\Big\} \to 0,
\eeq
which also implies, when $t \to 0$, that
\beq \label{xi}
\xi = \Phi^{-1}(1-t) \sim (\gamma \log(1/t))^{1/\gamma}.
\eeq

\subsection{Discovery times (LORD)} \label{sec:proof-LORD}

We apply LORD to the static setting under consideration.
Denote $\tau_l$ as the time of $l$-th discovery (with $\tau_0 = 0$), and $\Delta_l = \tau_l - \tau_{l-1}$ as the time between the $(l-1)$-th and $l$-th discoveries.  Assume a sequence satisfying \eqref{lambda} has been chosen.  Given the update rule of \eqref{LORD}, it can be seen that the inter-discovery times $\{\Delta_l : l \ge 1\}$ are IID.  

To prove \thmref{AGG-LORD}, we will use the following bound on the expected inter-discovery time.
 
\begin{prp} \label{prp:LORD}
Consider a static AGG mixture model with exponent $\gamma \ge 1$ parameterized as in \eqref{n-mu}. Assume that $\beta \in (0,1)$ and $r \ge 0$ are both fixed.  Assume that $r > \beta$ and let $\nu > 1$ be such that $\nu < r/\beta$.  If we apply LORD with $(\lambda_i)_{i = 1}^{\infty}$ defined as $\lambda_i \propto i^{-\nu}$ with $\sum_{i = 1}^{\infty} \lambda_i = q$, 
\beq
\E (\Delta_{l} \wedge n) \leq 2 n^{\beta} + C, \quad \text{for all } l >0,
\eeq
for some $C > 0$ that does not depend on $n$.  The same holds if we apply LORD with $(\lambda_i)_{i = 1}^{\infty}$ satisfying \eqref{newlambda} and $\sum_{i = 1}^{\infty} \lambda_i = q$.
\end{prp}

We prove this result. 
Recall the definition of $G$ in \eqref{G} and note that $G \ge \eps F$.
By the update rule of LORD algorithm, for all $m \geq 1$ we have 
\begin{align}
\P(\Delta_l > m) 
& = \prod_{i= \tau_{l-1} + 1}^{\tau_{l-1} + m} (1-G(\alpha_i)) 
= \prod_{i= \tau_{l-1} + 1}^{\tau_{l-1} + m} (1-G(\lambda_{i-\tau_{l-1}})) \\
&= \prod_{i=1}^{m} (1-G(\lambda_i)) 
\leq \exp \Big\{-\sum_{i=1}^{m} G(\lambda_i)\Big\} 
\leq \exp \Big\{-\eps \sum_{i=1}^{m} F(\lambda_i)\Big\}. 
\end{align}
Let $t^*$ be the value such that $\Phi^{-1}(1-t^*) = \mu$, i.e., $t^* = \Phi(-\mu) = n^{-r+o(1)}$ by the fact that $\Phi$ satisfies \defref{AGG}. Then, for $t \geq t^*$, we get
\beq
 \Phi^{-1}(1-t) \leq \Phi^{-1}(1-t^*) = \mu,
\eeq
and then
\beq
F(t) =\Phi(\mu - \Phi^{-1}(1-t)) \geq \Phi(\mu - \Phi^{-1}(1-t^*)) = \Phi(\mu - \mu) = \Phi(0) = 1/2,
\eeq
so that if $\lambda_i = Li^{-\nu} \geq t^*$, i.e., $i \leq n_1 := \lfloor (L/t^*)^{1/\nu} \rfloor = n^{r/\nu +o(1)}$, we have $F(\lambda_i) \geq \Phi(0) = 1/2$.

\begin{rem}\label{rem:lambda-proof1}
If instead $(\lambda_i)_{i = 1}^{\infty}$ satisfies \eqref{newlambda} then $i^{\nu} \lambda_i \to \infty$ as $i \to \infty$, so that exists a constant $L > 0$ such that $\lambda_i \ge L i^{-\nu}$ for all $i$, and this is all that we need to proceed.
\end{rem}

Thus, for $m \leq n_1$,
\beq
\sum_{i = 1}^{m} F(\lambda_i) \geq m/2,
\eeq
and for $m > n_1$,
\beq
\sum_{i = 1}^{m} F(\lambda_i) \geq \sum_{i = 1}^{n_1} F(\lambda_i) \geq n_1 /2.
\eeq
Thus,
\beq
\P(\Delta_l > m) \le \exp \{-\eps (m \wedge n_1)/2\}.
\eeq 

Next we bound $\E(\Delta_l \wedge n)$. 
Due to the fact that $\{\Delta_l \wedge n > m\} = \{\Delta_l > m\}$ for $1 \leq m \leq n-1$, and $\{\Delta_l \wedge n > m\} = \emptyset$ if $m \ge n$, we have 
\begin{align}
\E(\Delta_l \wedge n) 
&= \sum_{m = 0}^{\infty} \P (\Delta_l \wedge n > m) \\
&= \sum_{m = 1}^{n-1} \P (\Delta_l > m) + 1 \\
& \leq \sum_{m = 1}^{n-1} \exp \{-\eps (m \wedge n_1)/2\} + 1.
\end{align}
We split the summation over $1 \leq m \leq n_1$ and $n_1 +1 \leq m \leq n$ and derive the corresponding upper bound separately.  
For the first part, 
\begin{align}
\sum_{m = 1}^{n_1} \exp \{-\eps (m \wedge n_1)/2\} 
= \sum_{m = 1}^{n_1} \exp \{-\eps m/2\} 
\leq \frac{1}{\exp\{\eps/2\} -1 } < \frac{2}{\eps} = 2n^{\beta}.
\end{align}
For the second part, 
\begin{align}
\sum_{m = n_1 +1}^{n-1} \exp \{-\eps (m \wedge n_1)/2\} 
= \sum_{m = n_1 +1}^{n-1} \exp \{-\eps n_1/2\} 
\leq n \exp \{-\eps n_1/2\} = o(1),
\end{align}
since $\eps n_1 = n^{r/\nu - \beta + o(1)}$ and $\frac{r}{\nu} > \beta$.
Combining the above two bounds, we obtain 
\beq
\E(\Delta_l \wedge n)  \leq 2 n^\beta + o(1) + 1.
\eeq
This establishes \prpref{LORD}.


\subsection{Proof of \thmref{AGG-LORD}}\label{sec:proofLORD}
Note the number of false nulls is $m = |\F_n| = \eps n \sim n^{1-\beta}$. The false non-discovery rate of LORD (denoted $\fnr_n$) is as follows:
\begin{align}
\fnr_n & = \E \bigg(\frac{\sum_{i = 1}^{n} \IND {i \notin \cH_0(n): P_i \ge \alpha_i}}{m}\bigg) \\
& = \frac{\sum_{i = 1}^{n} \E [\E (\IND {i \notin \cH_0(n): P_i \ge \alpha_i} \mid \alpha_i)]}{m}\\
& = \frac{\sum_{i = 1}^{n} \E [\P (i \notin \cH_0(n) , P_i \ge \alpha_i \mid \alpha_i)]}{m}\\
& = \frac{\sum_{i = 1}^{n} \E [\eps \bar{F}(\alpha_i)]}{\eps n}\\
& = \frac{\sum_{i = 1}^{n} \E [\bar{F}(\alpha_i)]}{n}.
\end{align}
So it suffices to bound the RHS of the equation.

Let $D(n)$ be the number of discoveries in first $n$ hypotheses $\cH(n)$ by applying LORD with the sequence $(\lambda_i)$.
Let $\tilde{\Delta}_l = \Delta_l = \tau_l - \tau_{l-1}$, for $1 \le l \le  D(n)$, and $\tilde{\Delta}_{D(n)+1} = n-\tau_{D(n)}$. Due to the fact that $0 \le \tilde{\Delta}_l \leq (\Delta_l \wedge n)$, for $1 \le l \le D(n)+1$, we have  for any fixed $\delta >0$ ,
\beq
\P(\tilde{\Delta}_l \ge \frac{\E(\Delta_l \wedge n)}{\delta}) \le \frac{\E(\tilde{\Delta}_l)}{\E(\Delta_l \wedge n)} \cdot \delta \le \delta, \quad \text{for } 1 \le l \le  D(n)+1,
\eeq
by Markov Inequality.
Note that $(\Delta_l \wedge n)$'s are IID. We define $M_n := \left \lceil \E(\Delta)/\delta \right \rceil$, where $\Delta \overset{d}{=}\Delta_l \wedge n$ for all $l >0$.

For any $i \in \cH(n)$, there exists only one $j = j(i) \in \{1,2,\dots, D(n) + 1\}$ such that $i \in (\tau_{j-1}, \tau_j \wedge n]$, and 
\begin{align}
\E \big[\bar{F}(\alpha_i)\big] & = \E \big[\bar{F}(\alpha_i) \cdot \IND {\tilde{\Delta}_{j(i)} \ge M_n}\big] + \E \big[\bar{F}(\alpha_i) \cdot \IND {\tilde{\Delta}_{j(i)} < M_n} \big]\\
& \leq \delta + \E \big[\bar{F}(\alpha_i) \cdot \IND {\tilde{\Delta}_{j(i)} < M_n} \big],
\end{align}
so that 
\beq \label{upper}
\sum_{i = 1}^{n} \E \big[\bar{F}(\alpha_i)\big] \le n\delta +  \E \bigg[\sum_{i = 1}^{n}\bar{F}(\alpha_i) \cdot \IND {\tilde{\Delta}_{j(i)} < M_n}\bigg].
\eeq
By \prpref{LORD}, there is $C>0$ not depending on $n$ such that
\beq
\E(\Delta) \leq 2 n^{\beta} + C, \quad \text{for all } l >0.
\eeq
And thus, there is some $L'>0$ (constant in $n$) such that, for $1 \le i \le M_n$,
\beq\label{lambda-proof1}
\lambda_i = L i^{-\nu} \ge L \cdot (M_n)^{-\nu} = L \cdot \left \lceil \E(\Delta)/\delta \right \rceil  ^{-\nu} \ge L \cdot  \left \lceil (2n^{\beta}+C)/\delta \right \rceil  ^{-\nu} \ge L' n^{-\beta \nu}.
\eeq

\begin{rem}\label{rem:lambda-proof2}
If instead $(\lambda_i)_{i = 1}^{\infty}$ satisfies \eqref{newlambda} then $i^{\nu} \lambda_i \to \infty$ as $i \to \infty$, so that exists a constant $L > 0$ such that $\lambda_i \ge L i^{-\nu}$ for all $i$, and this is all that we need to proceed.
\end{rem}

Since $\bar{F}$ is a decreasing function, the second term in RHS of \eqref{upper} can be bounded as
\begin{align}
\E \bigg[\sum_{i = 1}^{n}\bar{F}(\alpha_i) \cdot \IND {\tilde{\Delta}_{j(i)} < M_n}\bigg] & = \E \bigg[\sum_{j = 1}^{D(n)+1} \sum_{i =\tau_{j-1} +1}^{\tau_j \wedge n}\bar{F}(\alpha_i) \cdot \IND {\tilde{\Delta}_{j} < M_n} \bigg]\\
& = \E \bigg[\sum_{j = 1}^{D(n)+1} \sum_{i =1}^{\tilde{\Delta}_j}\bar{F}(\lambda_i) \cdot \IND {\tilde{\Delta}_{j} < M_n} \bigg] \\
& \le \E \bigg[\sum_{j = 1}^{D(n)+1} \sum_{i =1}^{\tilde{\Delta}_j}\bar{F}(L'n^{-\beta \nu}) \cdot \IND {\tilde{\Delta}_{j} < M_n} \bigg] \\
& \le \E \bigg[\sum_{i = 1}^{n}\bar{F}(L'n^{-\beta \nu})\bigg] \le n \cdot \bar{F}(L'n^{-\beta \nu}).
\end{align}
Combining these bounds, we obtain
\beq
\fnr_n (\R) = \frac{\sum_{i = 1}^{n} \E [\bar{F}(\alpha_i)]}{n} \leq \delta + \bar{F}(L'n^{-\beta \nu}).
\eeq
Since $L' n^{-\beta \nu} \to 0$ as $n \to \infty$, by equation \eqref{xi} we have
\beq
\xi_n := \Phi^{-1}(1-L' n^{-\beta \nu}) = (\gamma \beta \nu \log n)^{1/\gamma}(1+o(1)), 
\eeq
so that
\begin{align}
\mu -\xi_n 
&= (\gamma r \log n)^{1/\gamma} - (\gamma \beta \nu \log n)^{1/\gamma}(1+o(1)) \\
&\sim (r^{1/\gamma} - (\beta \nu)^{1/\gamma}) (\gamma \log n)^{1/\gamma}
\to \infty, \quad \text{as } n \to \infty,
\end{align}
since $r > \beta \nu$.
Therefore, $\bar F(L'n^{-\beta \nu}) = \bar{\Phi}(\mu - \xi_n) \to 0$ as $n \to \infty$.
Hence, 
\begin{align}
\limsup_{n \to \infty} \fnr_n \le \delta.
\end{align}
This being true for any $\delta > 0$, necessarily, $\fnr_n \to 0$ as $n \to \infty$.
This establishes \thmref{AGG-LORD}.

\subsection{Discovery times (LOND)} \label{sec:proof-LOND}

We apply LOND to the static setting under consideration.
Denote $\tau_l$ as the time of $l$-th discovery (with $\tau_0 = 0$), and $\Delta_l = \tau_l - \tau_{l-1}$ as the time between the $(l-1)$-th and $l$-th discoveries.  Assume a sequence satisfying \eqref{lambda} has been chosen.  Given the update rule of \eqref{LOND}, it can be seen that the inter-discovery times $\{\Delta_l : l \ge 1\}$ are i.i.d..  

To prove \thmref{AGG-LOND}, we will use the following bound on the expected discovery times.
 
\begin{prp} \label{prp:LOND}
Consider a static AGG mixture model with exponent $\gamma \ge 1$ parameterized as in \eqref{n-mu}. Assume that $\beta \in (0,1)$ and $r \in [0,1]$ are both fixed.  For any $\nu > 1$, if we apply LOND with $(\lambda_i)_{i = 1}^{\infty}$ defined as $\lambda_i \propto i^{-\nu}$ with $\sum_{i = 1}^{\infty} \lambda_i = q$, 
	\beq
	\E (\tau_{l} \wedge n) \leq l \cdot n^{\beta + (\nu^{1/\gamma} - r^{1/\gamma})^\gamma + b_n}, \quad \text{for all } l >0,
	\eeq
where $b_n \to 0$ as $n \to \infty$.
\end{prp}

We now prove this result.
By the update rule of LOND algorithm, for all $l \ge 0$, and all $m \ge \tau_l + 1$, we have 
\beq
\P(\tau_{l+1} > m \mid \tau_l) = \prod_{i=\tau_l+1}^{m} (1-G((l+1)\lambda_i)) \leq \exp \{-\sum_{i =\tau_l+1}^{m} G((l+1)\lambda_i) \}.
\eeq
Note $\tau_l$ is the time of $l$-th discovery (with $\tau_0 = 0$) by LOND.
Let $\tilde{\tau}_l = \tau_l \wedge n$. If $\tilde{\tau}_l = n$, we have $\E (\tilde{\tau}_{l+1} \mid \tilde{\tau}_l) = n = \tilde{\tau}_l$.
Otherwise, if $\tilde{\tau}_l = \tau_l < n$,
\begin{align} 
\E (\tilde{\tau}_{l+1} \mid \tilde{\tau}_l) &= \tau_l + 1 + \sum_{m = \tau_l + 1}^{\infty}\P(\tau_{l+1} \wedge n > m \mid \tau_l) \\
&= \tau_l + 1 + \sum_{m = \tau_l + 1}^{n-1}\P(\tau_{l+1} > m \mid \tau_l) \\
& \le \tau_l + 1 + \sum_{m = \tau_l + 1}^{n}\exp \{-\sum_{i =\tau_l+1}^{m} G((l+1)\lambda_i)\} \\
&= \tilde{\tau}_l + 1 + \sum_{m = \tilde{\tau}_l + 1}^{n}\exp \{-\sum_{i =\tilde{\tau}_l+1}^{m} G((l+1)\lambda_i)\}. \label{bound}
\end{align}

Next we bound $\E (\tilde{\tau}_{l+1} \mid \tilde{\tau}_l)$. Let $t^*$ be the value such that $\Phi^{-1}(1-t^*) = \mu$, i.e., $t^* = \Phi(-\mu) = n^{-r+o(1)}$ by the fact that $\Phi$ satisfies \defref{AGG}. Then, for $t \geq t^*$, we get
\beq
\Phi^{-1}(1-t) \leq \Phi^{-1}(1-t^*) = \mu,
\eeq
and,
\beq
F(t) =\Phi(\mu - \Phi^{-1}(1-t)) \geq \Phi(\mu - \Phi^{-1}(1-t^*)) = \Phi(\mu - \mu) = \Phi(0) = 1/2,
\eeq
so that if $(l+1)\lambda_i = (l+1)Li^{-\nu} \ge t^*$, i.e., $i \leq n_1 := \lfloor ((l+1)L/t^*)^{1/\nu} \rfloor = n^{r/\nu +o(1)}$, we have $F((l+1)\lambda_i) \geq \Phi(0) = 1/2$.

We consider the following cases.
\paragraph{Case 1: $\tilde{\tau}_l<n_1 < n$.}  

In this case, for $\tilde{\tau}_l +1 \le m \le n_1$,
\beq
\sum_{i=\tilde{\tau}_l+1}^{m} G((l+1)\lambda_i) \geq \sum_{i = \tilde{\tau}_l+1}^{m} \eps F((l+1)\lambda_i) \geq \eps \cdot (m - \tilde{\tau}_l) /2,
\eeq
and for $m > n_1$,
\beq
\sum_{i=\tilde{\tau}_l+1}^{m} G((l+1)\lambda_i) \geq \sum_{i = \tilde{\tau}_l+1}^{m} \eps F((l+1)\lambda_i) \geq \sum_{i = \tilde{\tau}_l+ 1}^{m} \eps F((l+1)\lambda_m) = (m - \tilde{\tau}_l) \eps F((l+1)\lambda_m),
\eeq
since $F(x)$ is non-decreasing.

We split the summation in \eqref{bound} over $\tau_l + 1 \leq m \leq n_1$ and $n_1 +1 \leq m \leq n$ and derive the corresponding upper bound separately. For the first part,
\begin{align}
\sum_{m = \tilde{\tau}_l + 1}^{n_1}\exp \{-\sum_{i =\tilde{\tau}_l+1}^{m} G((l+1)\lambda_i)\} &\leq \sum_{m =\tilde{\tau}_l+1}^{n_1} \exp \{-\eps (m-\tilde{\tau}_l)/2\}
= \sum_{m =1}^{n_1 -\tilde{\tau}_l} \exp \{-\eps m/2\}\\
&\leq \frac{1}{\exp\{\eps/2\} -1 } < \frac{2}{\eps} = 2n^{\beta}.
\end{align}
For the second part,
\begin{align}
\sum_{m = n_1 +1}^{n}\exp \Big\{-\sum_{i =\tilde{\tau}_l+1}^{m} G((l+1)\lambda_i)\Big\}
&\le \sum_{m = n_1 +1}^{n} \exp \{-(m - \tilde{\tau}_l) \eps F((l+1)\lambda_m)\} \\
& \le \sum_{m = n_1 +1}^{n} \exp \{-(m - n_1) \eps F((l+1)\lambda_n)\} \\
& \le \sum_{m = 1}^{n -n_1} \exp \{-m \eps F((l+1)\lambda_n)\} \\
& \le  \frac{1}{\exp\{\eps F((l+1)\lambda_n)\} -1 } \\
& < \frac{1}{\eps F((l+1)\lambda_n)} \le \frac{1}{\eps F(\lambda_n)}.
\end{align}

\paragraph{Case 2: $n_1 \le \tilde{\tau}_l < n$.}  

For this case, we don't need to split the summation, since 
\begin{align}
\sum_{m = \tilde{\tau}_l + 1}^{n}\exp \{-\sum_{i =\tilde{\tau}_l+1}^{m} G((l+1)\lambda_i)\} &\le \sum_{m = \tilde{\tau}_l +1}^{n} \exp \{-(m - \tilde{\tau}_l) \eps F((l+1)\lambda_m)\} \\
& \le \sum_{m = 1}^{n -\tilde{\tau}_l} \exp \{-m \eps F((l+1)\lambda_n)\} \\
& < \frac{1}{\eps F((l+1)\lambda_n)} \le \frac{1}{\eps F(\lambda_n)}.
\end{align}

\paragraph{Case 3: $n_1 \ge n$.} 
Since $\tilde{\tau}_l <n \le n_1$, we have that 
\begin{align}
\sum_{m = \tilde{\tau}_l + 1}^{n}\exp \{-\sum_{i =\tilde{\tau}_l+1}^{m} G((l+1)\lambda_i)\} &\le \sum_{m = \tilde{\tau}_l +1}^{n} \exp \{-(m - \tilde{\tau}_l) \eps /2\} \\
& \le \sum_{m = 1}^{n -\tilde{\tau}_l} \exp \{-m \eps /2\} < \frac{2}{\eps} = 2n^{\beta}.
\end{align}
Combining all the cases, we obtain 
\beq
\E (\tilde{\tau}_{l+1} \mid \tilde{\tau}_l) \le \tilde{\tau}_l + 1 + \frac{2}{\eps} + \frac{1}{\eps F(\lambda_n)},
\eeq
where $F(\lambda_n) = \bar{\Phi}(\xi_n - \mu)$, and $\xi_n := \Phi^{-1}(1-\lambda_n)$.
Since $\lambda_n = L n^{-\nu} \to 0$ as $n \to 0$, by equation \eqref{xi}, we have $\xi_n \sim (\gamma \nu \log n)^{1/\gamma}$, so that
\begin{align}
\xi_n - \mu 
&= (\gamma \nu \log n)^{1/\gamma} (1+o(1))- (\gamma r \log  n)^{1/\gamma} \\
&\sim (\nu^{1/\gamma} - r^{1/\gamma}) (\gamma \log n)^{1/\gamma}
\to \infty, \quad \text{as } n \to \infty,
\end{align}
by the fact that $\nu > 1 \ge r$.
By \defref{AGG}, 
\begin{align}
F(\lambda_n) 
= \bar{\Phi}(\xi_n - \mu ) = \exp \Big\{-\frac{(\xi_n - \mu)^{\gamma}}{\gamma} (1+o(1))\Big\} 
= n^{-(\nu^{1/\gamma} - r^{1/\gamma})^\gamma + o(1)}.
\end{align}
Thus, when $n$ is large enough,
\beq
\E (\tilde{\tau}_{l+1} \mid \tilde{\tau}_l) \le \tilde{\tau}_l + 1 + \frac{2}{\eps} + \frac{1}{\eps F(\lambda_n)} \le \tilde{\tau}_l + n^{\beta + (\nu^{1/\gamma} - r^{1/\gamma})^\gamma + o(1)}, \quad \text{for all } l > 0,
\eeq
where the $o(1)$ is uniform in $l$, and this further implies that
\beq
\E (\tilde{\tau}_{l+1}) = \E[\E (\tilde{\tau}_{l+1} \mid \tilde{\tau}_l)] \le \E(\tilde{\tau}_l) + n^{\beta + (\nu^{1/\gamma} - r^{1/\gamma})^\gamma + o(1)}, \quad \text{for all } l > 0,
\eeq
so that 
\beq \label{discoverybound}
\E (\tau_{l} \wedge n) \leq l \cdot n^{\beta + (\nu^{1/\gamma} - r^{1/\gamma})^\gamma + o(1)}, \quad \text{for all } l > 0.
\eeq

\subsection{Proof of \thmref{AGG-LOND}} \label{sec:proof-AGG-LOND}
It suffices to consider the case where $r \in [0,1]$ since the observations from $\H_0$ almost never get substantially larger than $(\gamma \log n)^{1/\gamma}$.
For $r \in [0,1]$, if $r - (1 - r^{1/\gamma})^\gamma > \beta$, we can choose $\nu > 1$ close to 1 and $\eta > 0$ close to 0 such that $
r > \rho := \beta + (\nu^{1/\gamma} - r^{1/\gamma})^\gamma + \nu -1 + \eta$.
By \prpref{LOND}, when $n$ is large enough, 
\beq
\E (\tau_{l} \wedge n) \leq l \cdot n^{\beta + (\nu^{1/\gamma} - r^{1/\gamma})^\gamma + \eta}, \quad \text{for all } l > 0.
\eeq
Fix $\delta > 0$ and let $n_2:= \lceil n^{\beta + (\nu^{1/\gamma} - r^{1/\gamma})^\gamma + \eta}/\delta \rceil$. Note $n_2 = o(n)$, since $1\ge r>\beta + (\nu^{1/\gamma} - r^{1/\gamma})^\gamma + \eta$.

For $n_2 \le i \le n$, let $\zeta_i := i\, \delta\, n^{-\beta - (\nu^{1/\gamma} - r^{1/\gamma})^\gamma - \eta}$, we get
\begin{align}
\P (D(i) < \zeta_i) &= \P(\tau_{\lceil \zeta_i \rceil} > i) \le \P(\tau_{\lceil \zeta_i \rceil} \ge i) = \P(\tau_{\lceil \zeta_i \rceil} \wedge n \ge i)\\
&\le \frac{\E (\tau_{\lceil \zeta_i \rceil} \wedge n)}{i}
\le \frac{\lceil \zeta_i \rceil \cdot n^{\beta + (\nu^{1/\gamma} - r^{1/\gamma})^\gamma + \eta}}{i}\\
&< \frac{(\zeta_i+1) \cdot n^{\beta + (\nu^{1/\gamma} - r^{1/\gamma})^\gamma + \eta}}{i}\\
&= \delta + \frac{n^{\beta + (\nu^{1/\gamma} - r^{1/\gamma})^\gamma + \eta}}{i} <2\delta.
\end{align}
By Rule \eqref{LOND} defining the LOND algorithm, 
\beq
\E [\bar{F}(\alpha_i)] 
= \E [\bar{F}(\lambda_i(D(i-1) +1))]
\le \E [\bar{F}(\lambda_i D(i))],
\eeq
due to the fact that $D(i-1) +1 \ge D(i)$ and that $\bar{F}(x)$ is a non-increasing function, so that LOND's false non-discovery rate (denoted $\fnr_n$) is bounded as follows 
\begin{align}
\fnr_n = \frac1n \sum_{i = 1}^{n} \E [\bar{F}(\alpha_i)] 
\le \frac1n \sum_{i = 1}^{n} \E [\bar{F}(\lambda_i D(i))].
\end{align}
For $1\le i \le n_2$, 
\beq
\frac1n \sum_{i = 1}^{n_2} \E [\bar{F}(\lambda_i D(i))] \le \frac{n_2}{n}.
\eeq
And for $n_2 +1 \le i \le n$,
\begin{align}
\E [\bar{F}(\lambda_i D(i))] &= \E [\bar{F}(\lambda_i D(i)) \cdot \IND {D(i) < \zeta_i}] + \E [\bar{F}(\lambda_i D(i)) \cdot \IND {D(i) \ge \zeta_i}] \\
& \le 2\delta + \bar{F}(\lambda_i \zeta_i),
\end{align}
and since $\nu > 1$, we have
\beq
\lambda_i \zeta_i = L \delta \cdot i^{1-\nu} \cdot n^{-\beta - (\nu^{1/\gamma} - r^{1/\gamma})^\gamma - \eta} \ge L \delta \cdot n^{1-\nu} \cdot n^{-\beta - (\nu^{1/\gamma} - r^{1/\gamma})^\gamma - \eta} = \lambda_n \zeta_n,
\eeq
which implies that,
\beq
\frac1n \sum_{i = 1}^{n} \E [\bar{F}(\lambda_i D(i))] \le \frac{n_2}{n} +\frac{n-n_2}{n} (2\delta + \bar{F}(\lambda_n \zeta_n))
\le 2\delta + \bar{F}(\lambda_n \zeta_n) + o(1).
\eeq
Since $\lambda_n \zeta_n = L\delta n^{-\rho} \to 0$ as $n \to \infty$,  by equation \eqref{xi}
\beq
\xi_n := \Phi^{-1}(1-\lambda_n \zeta_n) = (\gamma \rho \log n)^{1/\gamma}(1+o(1)), 
\eeq
then
\begin{align}
\mu -\xi_n 
&= (\gamma r \log n)^{1/\gamma} - (\gamma \rho \log  n)^{1/\gamma}(1+o(1)) \\
&\sim (r^{1/\gamma} - \rho^{1/\gamma}) (\gamma \log n)^{1/\gamma}
\to \infty, \quad \text{as } n \to \infty,
\end{align}
since $r >  \rho$.
Therefore, $\bar F(\lambda_n \zeta_n) = \bar{\Phi}(\mu - \xi_n) \to 0$ as $n \to \infty$.
Hence, 
\begin{align}
\limsup_{n \to \infty} \fnr_n \le 2\delta.
\end{align}
This being true for any $\delta > 0$, necessarily, $\fnr_n \to 0$ as $n \to \infty$.

\section{Acknowledgement}
We would like to thank Jiaqi Gu from Department of Computer Science, University of California, Los Angeles, for his help with the large-sample size numerical experiments in \secref{numerics}.

\bibliographystyle{chicago}
\bibliography{ref}

\end{document}

%% file: notation.tex
\def\A{\mathscr{A}}
\def\F{\mathscr{F}}

\def\R{\mathscr{R}}

\def\H{\bbH}

\def\fdr{\textsc{fdr}}
\def\fdp{\textsc{fdp}}
\def\fnr{\textsc{fnr}}
\def\fnp{\textsc{fnp}}

\newtheorem{thm}{Theorem}
\newtheorem{prp}{Proposition}

\theoremstyle{remark}
\newtheorem{Def}{Definition}
\newtheorem{rem}{Remark}

\def\beq{\begin{equation}} 
\def\eeq{\end{equation}}
\def\beqn{\begin{eqnarray*}}
\def\eeqn{\end{eqnarray*}}
\def\Bitem{\begin{itemize}\setlength{\itemsep}{.2in}}
\def\bitem{\begin{itemize}\setlength{\itemsep}{.05in}}
\def\eitem{\end{itemize}}
\def\Benum{\begin{enumerate}\setlength{\itemsep}{.2in}}
\def\benum{\begin{enumerate}\setlength{\itemsep}{.05in}}
\def\eenum{\end{enumerate}}
\def\bmult{\begin{multline*}}
\def\emult{\end{multline*}}
\def\bcenter{\begin{center}}
\def\ecenter{\end{center}}
\def\bframe{\begin{frame}}
\def\eframe{\end{frame}}

\newcommand{\thmref}[1]{Theorem~\ref{thm:#1}}
\newcommand{\prpref}[1]{Proposition~\ref{prp:#1}}

\newcommand{\secref}[1]{Section~\ref{sec:#1}}
\newcommand{\figref}[1]{Figure~\ref{fig:#1}}

\newcommand{\remref}[1]{Remark~\ref{rem:#1}}
\newcommand{\defref}[1]{Definition~\ref{def:#1}}

\def\cH{\mathcal{H}}

\def\bX{\mathbf{X}}

\def\bbH{\mathbb{H}}
\def\bbI{\mathbb{I}}

\newcommand{\E}{\operatorname{\mathbb{E}}}
\renewcommand{\P}{\operatorname{\mathbb{P}}}

\def\eps{\varepsilon}

\def\1{\mathbbm{1}}
\newcommand{\IND}[1]{\bbI\{ #1 \}}

%% file: sequential-multiple-testing.bbl
\begin{thebibliography}{}

\bibitem[\protect\citeauthoryear{Aharoni and Rosset}{Aharoni and
  Rosset}{2014}]{aharoni2014generalized}
Aharoni, E. and S.~Rosset (2014).
\newblock Generalized $\alpha$-investing: definitions, optimality results and
  application to public databases.
\newblock {\em Journal of the Royal Statistical Society: Series B (Statistical
  Methodology)\/}~{\em 76\/}(4), 771--794.

\bibitem[\protect\citeauthoryear{Arias-Castro and Chen}{Arias-Castro and
  Chen}{2016}]{ariaschen2016distribution}
Arias-Castro, E. and S.~Chen (2016).
\newblock Distribution-free multiple testing.
\newblock {\em arXiv preprint arXiv:1604.07520\/}.

\bibitem[\protect\citeauthoryear{Bartroff}{Bartroff}{2014}]{bartroff2014multiple}
Bartroff, J. (2014).
\newblock Multiple hypothesis tests controlling generalized error rates for
  sequential data.
\newblock {\em arXiv preprint arXiv:1406.5933\/}.

\bibitem[\protect\citeauthoryear{Bartroff and Song}{Bartroff and
  Song}{2013}]{bartroff2013sequential}
Bartroff, J. and J.~Song (2013).
\newblock Sequential tests of multiple hypotheses controlling false discovery
  and nondiscovery rates.
\newblock {\em arXiv preprint arXiv:1311.3350\/}.

\bibitem[\protect\citeauthoryear{Bartroff and Song}{Bartroff and
  Song}{2014}]{bartroff2014sequential}
Bartroff, J. and J.~Song (2014).
\newblock Sequential tests of multiple hypotheses controlling type i and ii
  familywise error rates.
\newblock {\em Journal of statistical planning and inference\/}~{\em 153},
  100--114.

\bibitem[\protect\citeauthoryear{Benjamini and Hochberg}{Benjamini and
  Hochberg}{1995}]{benjamini1995controlling}
Benjamini, Y. and Y.~Hochberg (1995).
\newblock Controlling the false discovery rate: A practical and powerful
  approach to multiple testing.
\newblock {\em Journal of the Royal Statistical Society. Series B
  (Methodological)\/}~{\em 57\/}(1), 289--300.

\bibitem[\protect\citeauthoryear{Bogdan, Chakrabarti, Frommlet, and
  Ghosh}{Bogdan et~al.}{2011}]{bogdan2011asymptotic}
Bogdan, M., A.~Chakrabarti, F.~Frommlet, and J.~K. Ghosh (2011).
\newblock Asymptotic bayes-optimality under sparsity of some multiple testing
  procedures.
\newblock {\em The Annals of Statistics\/}, 1551--1579.

\bibitem[\protect\citeauthoryear{Butucea, Stepanova, and Tsybakov}{Butucea
  et~al.}{2015}]{butucea2015variable}
Butucea, C., N.~A. Stepanova, and A.~B. Tsybakov (2015).
\newblock Variable selection with hamming loss.
\newblock {\em arXiv preprint arXiv:1512.01832\/}.

\bibitem[\protect\citeauthoryear{Dickhaus}{Dickhaus}{2014}]{dickhaus2014simultaneous}
Dickhaus, T. (2014).
\newblock {\em Simultaneous statistical inference}.
\newblock Springer.

\bibitem[\protect\citeauthoryear{Donoho and Jin}{Donoho and
  Jin}{2004}]{donoho2004higher}
Donoho, D. and J.~Jin (2004).
\newblock Higher criticism for detecting sparse heterogeneous mixtures.
\newblock {\em The Annals of Statistics\/}~{\em 32\/}(3), 962--994.

\bibitem[\protect\citeauthoryear{Dudoit and van~der Laan}{Dudoit and van~der
  Laan}{2007}]{dudoit2007multiple}
Dudoit, S. and M.~J. van~der Laan (2007).
\newblock {\em Multiple testing procedures with applications to genomics}.
\newblock Springer Science \& Business Media.

\bibitem[\protect\citeauthoryear{Fithian, Sun, and Taylor}{Fithian
  et~al.}{2014}]{fithian2014optimal}
Fithian, W., D.~Sun, and J.~Taylor (2014).
\newblock Optimal inference after model selection.
\newblock {\em arXiv preprint arXiv:1410.2597\/}.

\bibitem[\protect\citeauthoryear{Fithian, Taylor, Tibshirani, and
  Tibshirani}{Fithian et~al.}{2015}]{fithian2015selective}
Fithian, W., J.~Taylor, R.~Tibshirani, and R.~Tibshirani (2015).
\newblock Selective sequential model selection.
\newblock {\em arXiv preprint arXiv:1512.02565\/}.

\bibitem[\protect\citeauthoryear{Foster and Stine}{Foster and
  Stine}{2008}]{foster2008alpha}
Foster, D.~P. and R.~A. Stine (2008).
\newblock $\alpha$-investing: a procedure for sequential control of expected
  false discoveries.
\newblock {\em Journal of the Royal Statistical Society: Series B (Statistical
  Methodology)\/}~{\em 70\/}(2), 429--444.

\bibitem[\protect\citeauthoryear{Foygel-Barber and Cand{\`e}s}{Foygel-Barber
  and Cand{\`e}s}{2015}]{barber2015controlling}
Foygel-Barber, R. and E.~J. Cand{\`e}s (2015).
\newblock Controlling the false discovery rate via knockoffs.
\newblock {\em The Annals of Statistics\/}~{\em 43\/}(5), 2055--2085.

\bibitem[\protect\citeauthoryear{Genovese and Wasserman}{Genovese and
  Wasserman}{2002}]{genovese2002operating}
Genovese, C. and L.~Wasserman (2002).
\newblock Operating characteristics and extensions of the false discovery rate
  procedure.
\newblock {\em Journal of the Royal Statistical Society: Series B (Statistical
  Methodology)\/}~{\em 64\/}(3), 499--517.

\bibitem[\protect\citeauthoryear{G'Sell, Wager, Chouldechova, and
  Tibshirani}{G'Sell et~al.}{2016}]{g2016sequential}
G'Sell, M.~G., S.~Wager, A.~Chouldechova, and R.~Tibshirani (2016).
\newblock Sequential selection procedures and false discovery rate control.
\newblock {\em Journal of the Royal Statistical Society: Series B (Statistical
  Methodology)\/}~{\em 78\/}(2), 423--444.

\bibitem[\protect\citeauthoryear{Ingster}{Ingster}{1997}]{ingster1997some}
Ingster, Y.~I. (1997).
\newblock Some problems of hypothesis testing leading to infinitely divisible
  distributions.
\newblock {\em Mathematical Methods of Statistics\/}~{\em 6\/}(1), 47--69.

\bibitem[\protect\citeauthoryear{Ingster and Suslina}{Ingster and
  Suslina}{2003}]{IngsterBook}
Ingster, Y.~I. and I.~A. Suslina (2003).
\newblock {\em Nonparametric goodness-of-fit testing under {G}aussian models},
  Volume 169 of {\em Lecture Notes in Statistics}.
\newblock New York: Springer-Verlag.

\bibitem[\protect\citeauthoryear{Javanmard and Montanari}{Javanmard and
  Montanari}{2015}]{javanmard2015online}
Javanmard, A. and A.~Montanari (2015).
\newblock On online control of false discovery rate.
\newblock {\em arXiv preprint arXiv:1502.06197\/}.

\bibitem[\protect\citeauthoryear{Javanmard and Montanari}{Javanmard and
  Montanari}{2016}]{javanmard2016online}
Javanmard, A. and A.~Montanari (2016).
\newblock Online rules for control of false discovery rate and false discovery
  exceedance.
\newblock {\em arXiv preprint arXiv:1603.09000\/}.

\bibitem[\protect\citeauthoryear{Ji, Jin, et~al.}{Ji et~al.}{2012}]{ji2012ups}
Ji, P., J.~Jin, et~al. (2012).
\newblock Ups delivers optimal phase diagram in high-dimensional variable
  selection.
\newblock {\em The Annals of Statistics\/}~{\em 40\/}(1), 73--103.

\bibitem[\protect\citeauthoryear{Jin and Ke}{Jin and Ke}{2014}]{jin2014rare}
Jin, J. and T.~Ke (2014).
\newblock Rare and weak effects in large-scale inference: methods and phase
  diagrams.
\newblock {\em arXiv preprint arXiv:1410.4578\/}.

\bibitem[\protect\citeauthoryear{Lei and Fithian}{Lei and
  Fithian}{2016}]{lei2016power}
Lei, L. and W.~Fithian (2016).
\newblock Power of ordered hypothesis testing.
\newblock {\em arXiv preprint arXiv:1606.01969\/}.

\bibitem[\protect\citeauthoryear{Li and Barber}{Li and
  Barber}{2016}]{li2016accumulation}
Li, A. and R.~F. Barber (2016).
\newblock Accumulation tests for fdr control in ordered hypothesis testing.
\newblock {\em Journal of the American Statistical
  Association\/}~(just-accepted), 1--38.

\bibitem[\protect\citeauthoryear{Lockhart, Taylor, Tibshirani, and
  Tibshirani}{Lockhart et~al.}{2014}]{lockhart2014significance}
Lockhart, R., J.~Taylor, R.~J. Tibshirani, and R.~Tibshirani (2014).
\newblock A significance test for the lasso.
\newblock {\em Annals of statistics\/}~{\em 42\/}(2), 413.

\bibitem[\protect\citeauthoryear{Meinshausen, Maathuis, B{\"u}hlmann,
  et~al.}{Meinshausen et~al.}{2011}]{meinshausen2011asymptotic}
Meinshausen, N., M.~H. Maathuis, P.~B{\"u}hlmann, et~al. (2011).
\newblock Asymptotic optimality of the westfall--young permutation procedure
  for multiple testing under dependence.
\newblock {\em The Annals of Statistics\/}~{\em 39\/}(6), 3369--3391.

\bibitem[\protect\citeauthoryear{Neuvial and Roquain}{Neuvial and
  Roquain}{2012}]{neuvial2012false}
Neuvial, P. and E.~Roquain (2012).
\newblock On false discovery rate thresholding for classification under
  sparsity.
\newblock {\em The Annals of Statistics\/}~{\em 40\/}(5), 2572--2600.

\bibitem[\protect\citeauthoryear{Roquain}{Roquain}{2011}]{roquain2011type}
Roquain, E. (2011).
\newblock Type i error rate control in multiple testing: a survey with proofs.
\newblock {\em Journal de la Soci\'et\'e Fran\c{c}aise de Statistique\/}~{\em
  152\/}(2), 3--38.

\bibitem[\protect\citeauthoryear{Storey}{Storey}{2007}]{storey2007optimal}
Storey, J.~D. (2007).
\newblock The optimal discovery procedure: a new approach to simultaneous
  significance testing.
\newblock {\em Journal of the Royal Statistical Society: Series B (Statistical
  Methodology)\/}~{\em 69\/}(3), 347--368.

\bibitem[\protect\citeauthoryear{Sun and Cai}{Sun and
  Cai}{2007}]{sun2007oracle}
Sun, W. and T.~T. Cai (2007).
\newblock Oracle and adaptive compound decision rules for false discovery rate
  control.
\newblock {\em Journal of the American Statistical Association\/}~{\em
  102\/}(479), 901--912.

\end{thebibliography}
